\documentclass[11pt]{amsart}
\theoremstyle{plain}
\newtheorem{thm}{Theorem}[section]
\newtheorem{theorem}[thm]{Theorem}

\newtheorem{lemma}[thm]{Lemma}
\newtheorem{corollary}[thm]{Corollary}
\newtheorem{proposition}[thm]{Proposition}
\theoremstyle{definition}
\newtheorem{remark}[thm]{Remark}

\newtheorem{notation}[thm]{Notation}

\newtheorem{definition}[thm]{Definition}

\newtheorem{example}[thm]{Example}

\newtheorem{question}[thm]{Question}

\numberwithin{equation}{section}
\newcommand{\sC}{{\mathcal C}}
\newcommand{\sD}{{\mathcal D}}
\newcommand{\sE}{{\mathcal E}}
\newcommand{\sF}{{\mathcal F}}

\newcommand{\sH}{{\mathcal H}}

\newcommand{\sJ}{{\mathcal J}}
\newcommand{\sK}{{\mathcal K}}

\newcommand{\sO}{{\mathcal O}}

\newcommand{\sT}{{\mathcal T}}
\newcommand{\sU}{{\mathcal U}}
\newcommand{\sV}{{\mathcal V}}

\newcommand{\A}{{\mathbb A}}

\newcommand{\C}{{\mathbb C}}

\newcommand{\BP}{{\mathbb P}}

\newcommand{\Q}{{\mathbb Q}}

\newcommand\zero{\rm Zero}

\def\Sym{\mathop{\rm Sym}\nolimits}

\def\Hom{\mathop{\rm Hom}\nolimits}

\title[Unbendable rational curves]{Unbendable rational curves of Goursat type and Cartan type}

\author{Jun-Muk Hwang and Qifeng Li}

\thanks{This work was supported by the Institute for Basic Science (IBS-R032-D1).}

\begin{document}

\maketitle

\begin{abstract}
We study unbendable rational curves, i.e., nonsingular rational curves in a complex manifold of dimension $n$ with normal bundles isomorphic to $$\sO_{\BP^1}(1)^{\oplus p} \oplus  \sO_{\BP^1}^{\oplus (n-1-p)}$$ for some nonnegative integer $p$. Well-known examples arise from algebraic geometry as general minimal rational curves of  uniruled projective manifolds.
After describing the relations between the differential geometric properties of the natural distributions on the deformation spaces of unbendable rational curves and the projective  geometric properties of their varieties of minimal rational tangents, we concentrate on the case of $p=1$ and $n \leq 5$, which is the simplest nontrivial situation.  In this case, the families of unbendable rational curves  fall  essentially into  two classes: Goursat type or Cartan type. Those of Goursat type arise from ordinary differential equations and those of Cartan type have special features related to contact geometry.
We show that the family of lines on any nonsingular cubic 4-fold is of Goursat type, whereas the family of lines on
a general quartic 5-fold is of Cartan type, in the proof of which the projective geometry of varieties of minimal rational tangents plays a key role.
\end{abstract}

\medskip
MSC2010: 58A30, 32C25, 14J70

\section{Introduction}
A nonsingular rational curve in a complex manifold of dimension $n$ is said to be unbendable if its normal bundle is
isomorphic to $$\sO_{\BP^1}(1)^{\oplus p} \oplus  \sO_{\BP^1}^{\oplus (n-1-p)}$$ for some nonnegative integer $p$.
There are many examples of
unbendable rational curves arising from algebraic geometry. Any uniruled projective manifold  has minimal rational curves and general minimal rational curves are unbendable. In this case, unbendable rational curves are sometimes called standard rational curves, indicating that they are general among minimal rational curves (see Section 1.1 of  \cite{HM} or Section 1.2 of \cite{Hw01}). Understanding the germ of an unbendable rational curve in a complex manifold, i.e. the biholomorphic structure of neighborhoods of such a curve, is important especially when they are minimal rational curves in a Fano manifold of Picard number 1, because quite often the germ of the curve determines the biregular type of  the ambient Fano manifold by Cartan-Fubini type extension theorem (see Section 3 of \cite{HM} or Section 3 of \cite{Hw01}).

One approach to study the germ of an unbendable rational curve is by looking at certain distributions (i.e. Pfaffian  systems) on
the deformation space of the unbendable rational curves in the ambient complex manifold (i.e. the corresponding open subset in the Douady space of the ambient manifold). In many cases, the germ of an unbendable rational curve is determined by the germ of these distributions in a neighborhood of the corresponding point in the Douady space. In this paper, we employ this correspondence to study the germ of unbendable rational curves.

The relation between the geometry of rational curves on a complex manifold and certain natural differential systems on the corresponding Douady space has been  studied much in twistor theory (e.g. \cite{Hi}).
What is new in our approach is the role of the varieties of minimal rational tangents (VMRT in abbreviation), the complex submanifold in the projectivized tangent bundle of the ambient complex manifold traced by tangents to the deformations of the unbendable rational curves. The original notion of VMRT (e.g. in \cite{HM}) is for minimal rational curves on uniruled projective manifolds and they are defined as certain projective subvarieties in the projectivized tangent spaces of the uniruled projective manifolds. In this paper, we consider   VMRT for unbendable rational curves and they are defined (see Definition \ref{d.vmrt}) as certain complex submanifolds  in the projectivized tangent spaces. In fact, for our purpose, it is sufficient to consider the germs of VMRT. We show that   the projective geometry of VMRT is intricately related to the properties of the distributions on the deformation space, as explained in Section \ref{s.vmrt}.

To explain our main results, let us recall briefly  the notion of the growth vector of a distribution $D$ on a complex manifold $M$ (see Section \ref{ss.d} for a precise definition).
Associated to $D \subset TM$ is a sequence of saturated subsheaves $$D = \partial^0 D \subset \partial^1 D \subset \partial^2 D \subset \cdots \subset \partial^k D = \partial^{k+1} D$$ of $TM$ for some nonnegative integer $k$,  generated by the successive brackets of local vector fields belonging to $D$. The growth vector of $D$ is the strictly increasing sequence of positive integers $$( {\rm rank} (D), {\rm rank}(\partial^1 D), \ldots, {\rm rank}(\partial^k D)).$$  It is the most basic invariant of a distribution. We say that $D$ is bracket-generating if $\partial^k D = TM$, i.e.,  the last entry of its  growth vector is $\dim M$.

In the current article, we concentrate on
 the case $p =1$, i.e. when there is exactly one $\sO(1)$-factor in the normal bundle of  unbendable rational curves. This is the case when the anti-canonical degree of the rational curves is 3.
In this case, the deformation space has a natural rank 2 distribution and the germ of unbendable rational curves can be recovered from some extra structure on the germ of this rank 2 distribution (Proposition \ref{p.rank2}). There are two particularly well-understood classes of rank 2 distributions:  rank 2 distributions $D \subset TM$ satisfying $${\rm rank}(\partial^{i}D) = {\rm rank}(\partial^{i-1}  D) + 1 \mbox{ for all } 1 \leq i \leq k \mbox{ and }  \partial^k D = TM,$$ which are called Goursat distributions (see \cite{MZ}) and  distributions with the growth vector $(2,3,5), $ which \'{E}. Cartan studied extensively in  \cite{Ca}.  We investigate the geometry of unbendable rational curves whose associated distributions belong to one of these two classes of distributions. Let us call them unbendable rational curves of  Goursat type and of Cartan type, respectively.

Unbendable rational curves of Goursat type arise from ordinary differential equations and we describe precisely which ordinary differential equations give rise to unbendable rational curves (see Theorem \ref{t.ode} for a precise version):

\begin{theorem}\label{t.1Goursat}
Unbendable rational curves of Goursat type arise from ordinary differential equations of the type
$$\frac{{\rm d}^n u}{{\rm d} t^n}  = a_3 (\frac{{\rm d}^{n-1} u}{{\rm d} t^{n-1}})^3 + a_2 (\frac{{\rm d}^{n-1} u}{{\rm d} t^{n-1}})^2 + a_1 \frac{{\rm d}^{n-1} u}{{\rm d} t^{n-1}} + a_0$$ where    $a_0, a_1, a_2, a_3$ are local holomorphic functions of the $n$ variables  $t$, $u$, $\frac{{\rm d} u}{{\rm d} t}$, $\ldots,$ $\frac{{\rm d}^{n-2} u}{{\rm d} t^{n-2}}$. \end{theorem}

When $n=2$, Theorem \ref{t.1Goursat} is precisely Theorem 3.1 of \cite{Hi}.

To study unbendable rational curves of Cartan type, we use some structure theory of $(2,3,5)$-distributions (e.g. \cite{BH}, \cite{Z99}, \cite{Z06}). An interesting consequence is the following   result
 (see Theorem \ref{t.Zelenko} and Theorem \ref{t.model} for a precise version):

\begin{theorem}\label{t.Z}
There is a natural 1-1 correspondence between  germs of a $(2,3,5)$-distribution at  general points  and germs of general members of a bracket-generating family of unbendable rational curves of anti-canonical degree 3 in a complex manifold of dimension $5$ equipped with a contact structure such that the rational curves are tangent to the contact structure. \end{theorem}

Here, a family of unbendable rational curves is bracket-generating if the associated distribution in the deformation space of the curves is a bracket-generating distribution.
There is a well-known  geometric interpretation of  $(2,3,5)$-distributions  as the mechanical systems of rolling two surfaces without slipping or twisting in the three-dimensional Euclidean space (see Section 4.4 of \cite{BH} or Section 6.8 of \cite{Mo}). Theorem \ref{t.Z} can be viewed as another geometric interpretation of  $(2,3,5)$-distributions in the context of complex geometry.

For unbendable rational curves arising from algebraic geometry as minimal rational curves, it is often not easy to determine the differential geometric properties of the associated distributions. Even the growth vectors of the distributions, the simplest invariants of distributions, are not easy to compute. The simplest nontrivial question in this direction is to determine whether a given family of unbendable rational curves with $p=1$ in a complex manifold of dimension 5 is of Goursat type or of Cartan type. Our key result is the following criterion (see Theorem \ref{t.test5} and Theorem \ref{t.test}):

\begin{theorem}\label{t.testA} A bracket-generating family of unbendable rational curves with $p=1$ in a 5-dimensional manifold is of Goursat type if and only if the third fundamental forms of the VMRT along a general member is zero at some point. If it is not of Goursat type, then it is of Cartan type. \end{theorem}

Using this criterion, we prove the following (see Theorem \ref{t.general}).

\begin{theorem}\label{t.line}
The family of lines on a general hypersurface of degree 4 in $\BP^6$ is a family of unbendable rational curves of Cartan type.
\end{theorem}

This illustrates the advantage of using VMRT in our approach. We should mention, however, that the computation can be still tricky even after using Theorem \ref{t.testA}. For the proof of Theorem \ref{t.line}, we need to select a  hypersurface with
a special involution to simplify the computation.

 In this context,
Theorem \ref{t.line} leads to the following two questions which we leave for future studies.

\begin{question} Is the family of lines on {\em any} nonsingular hypersurface of degree $4$ in $\BP^6$ of Cartan type? \end{question}

\begin{question} What is the growth vector of the rank 2 distribution $\sD$ in Definition \ref{d.sD} for the family of lines on a general hypersurface of degree $n-1$ in $\BP^{n+1}$? \end{question}

 This paper is organized as follows. In Section \ref{s.prelim}, we collect standard facts on some differential geometric notions to fix our terminology and notation. In Section \ref{s.vmrt}, we present a general theory of the natural distribution on the deformation space of unbendable rational curves and its relation to VMRT. Section \ref{s.Goursat} studies unbendable rational curves of Goursat type and Section \ref{s.Cartan} studies unbendable rational curves of Cartan type. Finally, in Section \ref{s.quartic}, we prove Theorem \ref{t.line}.

\section{Preliminaries}\label{s.prelim}

We work in the complex analytic category: all  manifolds, bundles  and maps are holomorphic. Open subsets are taken in  Euclidean topology, unless stated otherwise. The tangent bundle of a complex manifold $M$ is denoted by $TM$. For a vector bundle $W$ of rank $r$ on $M$, its projectivization $\BP W \to M$ is the $\BP^{r-1}$-bundle whose fibers consist of 1-dimensional subspaces in  fibers of $W \to M$. The tautological line bundle on $\BP W$ is denoted by $\sO_{\BP W}(-1)$ and its dual bundle is denoted by $\sO_{\BP W}(1).$

\subsection{Distributions}\label{ss.d}

A {\em distribution} $D$ on a complex manifold $M$ is a vector subbundle $D \subset TM$ of the tangent bundle. By abuse of notation, we denote by $D \subset TM$ also the locally free subsheaf of local sections of $D$, a subsheaf of the sheaf of vector fields on $M$.
A {\em meromorphic distribution} $D$ of rank $r$ on a complex manifold $M$ is a coherent subsheaf $D$  of rank $r$ in the sheaf of vector fields on $M$ which is saturated in the sense that there exists no subsheaf $D'$ of the sheaf of vector fields satisfying $D \subset D', D \neq D'$ and ${\rm rank} (D) = {\rm rank} (D')$. By abuse of notation, we write $D \subset TM$ to denote a meromorphic distribution. When the restriction of $D$ to a neighborhood of a point $x \in M$ is a distribution, we write $D_x \subset T_xM$ to  denote the fiber of the corresponding vector subbundle at $x$.  The saturation condition implies that a meromorphic distribution induces a distribution on a Zariski open subset of $M$ which is the complement of a closed analytic subset of codimension at least 2. It implies also that $D$   is determined by the fibers $D_x \subset T_xM$ at general points $x \in M$.    Given a  distribution $D \subset TM$, Lie brackets of local sections $[,]: \wedge^2 D \to TM$ induces   a  homomorphism of vector bundles $${\rm Levi}^D: \wedge^2 D \to TM/D,$$ called  the {\em Levi tensor} of $D$.

The {\em Cauchy characteristic} of a meromorphic distribution $D$ on $M$ is
the subsheaf ${\rm Ch}(D) \subset D$ defined as follows: a local section $v$ of $D$ belongs to ${\rm Ch}(D)$ if and only if  $[v, w]$ belongs to $D$ for any local section $w$ of $D$. The Cauchy characteristic is involutive, defining a meromorphic foliation.

Associated to a meromorphic distribution $D \subset TM$ are the derived meromorphic distributions $\partial^i D$ and $\partial^{(i)} D$ defined inductively by
$$ \partial^1 D = \partial^{(1)} D = [D,D] + D,$$ $$\partial^{i+1} D = [\partial^{i} D, \partial^i D] + \partial^i D, \mbox{ and }
\partial^{(i+1)} D = [D, \partial^{(i)} D] + \partial^{(i)} D.$$ Writing $\partial^{-1} D =0$ and $ \partial^0 D = D$, we have an integer $d \geq 0$ such that  $${\rm rank}(\partial^{d-1}D) \neq {\rm rank}(\partial^{d} D) = {\rm rank}( \partial^{d+i} D)$$ for any $i >0$. The strictly increasing sequence of integers $$({\rm rank}(D), {\rm rank}(\partial D), {\rm rank}(\partial^2 D), \ldots, {\rm rank} (\partial^d D))$$ is called the {\em growth vector} of $D$.  If ${\rm rank}(\partial^d D)= \dim M$, we say that $D$ is {\em bracket-generating}.

We say that a meromorphic distribution $D \subset TM$ is {\em regular} at a point $x \in M$ (equivalently,  $x$ is a {\em regular point} of $D$) if
 ${\rm Ch}(D), \partial^i D$ and $\partial^{(i)} D$ for all $i \geq 0$ induce distributions in a neighborhood  of $x \in M$.
 Regular points of a meromorphic distribution  form a Zariski-open subset in $M$ called the {\em regular locus} of $D$.
 We say that $D$ is a {\em regular distribution} if every point of $M$ is a regular point of $D$.

A vector field $v$ on $M$ is an {\em infinitesimal automorphism} of a meromorphic distribution $D$, if the  1-parameter germs of local biholomorphisms of $M$ (and consequently $TM$) generated by $v$ preserve $D.$  The Lie algebra of all infinitesimal automorphisms of $D$ is denoted by $\frak{a}\frak{u}\frak{t}(D)$.

When $\pi: Y \to M$ is a submersion of complex manifolds, denote by $T^{\pi} \subset TY$ the distribution given by ${\rm Ker}({\rm d} \pi)$.
For a meromorphic distribution $D$ on $M$,
 the {\em inverse image} $\pi^{-1}D$ is a meromorphic distribution on $Y$ whose fiber at a general point
$y \in Y$ is $$ (\pi^{-1}D)_y = ({\rm d}_y \pi)^{-1}(D_{\pi(y)})$$ where ${\rm d}_y \pi : T_y Y \to T_{\pi(y)} M$ is the differential of $\pi$ at $y$.  It is easy to see that  \begin{equation}\label{e.dist}  \partial^{k} (\pi^{-1}D) = \pi^{-1} (\partial^k D), \ \partial^{(k)}( \pi^{-1}D) = \pi^{-1} (\partial^{(k)} D)\end{equation} for any $k \geq 1$ and ${\rm Ch}(\pi^{-1}D) = \pi^{-1}{\rm Ch}(D)$.

For a distribution $D \subset TM$, let $\pi: \BP D \to M$ be the projectivization of the vector bundle $D$. Define the distribution ${\rm pr}(D) $ on $\BP D$,  called the {\em prolongation } of $D$, whose fiber at $[v] \in \BP D$ corresponding to $0 \neq v \in D$ is $$ {\rm pr}(D)_{[v]} := ({\rm d}_{[v]} \pi)^{-1} (\C v).$$ Then \begin{equation}\label{e.pr} \partial ({\rm pr}(D)) = \pi^{-1}D \end{equation} from Proposition 5.1 of \cite{MZ} (this is also a special case of Proposition 1 in \cite{HM04}).

\subsection{Fundamental forms of projective submanifolds}\label{ss.p}
A good reference for fundamental forms is Section 12.1 of \cite{IL}, where they are formulated in terms of differential forms. For our purpose, it is convenient to reformulate them in terms of vector fields as follows.

For a vector space $V$, its projectivization $\BP V$ is the set of 1-dimensional subspaces of $V$.
For a (not necessarily closed) complex submanifold $Z \subset \BP V$, denote by $\hat{Z} \subset V$ be the  union of 1-dimensional subspaces corresponding to $Z$.
The affine tangent space $\hat{T}_z Z \subset V$ at a point $z \in Z$ is the tangent space of $\hat{Z}$ at a point of $\hat{z} \setminus 0$. There is a natural identification
$$T_z Z = \Hom (\hat{z}, \hat{T}_z Z/\hat{z}), \ \mbox{ (equivalently, } \hat{T}_z Z/\hat{z} = \hat{z} \otimes T_z Z) ).$$
Let $\xi: V\setminus 0\rightarrow\BP V$ be the natural $\C^*$-bundle. Fix a  linear coordinate system $x_1, \ldots, x_n$ on $V$ where $n= \dim V$.  Let $v_1,\ldots, v_k$ be local vector fields on $Z$ near a point $z \in Z.$ Given  any point $z' \in V \setminus 0$ with $\xi(z') = z$, we can choose holomorphic functions $v_{i,j} $ in a neighborhood of $z'$ in $V$ such that  $$\hat{v}_i=\sum_{j=1}^n v_{i, j}\frac{\partial}{\partial x_j}$$ are local vector fields tangent to $\hat{Z}$ near $z'$  satisfying  ${\rm d} \xi (\hat{v}_i )|_{\hat{Z}} = v_i$.  Write $$ \hat{v}_2(\hat{v}_1) = \sum_{j=1}^n\hat{v}_2(v_{1, j})\frac{\partial}{\partial x_j}$$ and inductively $$ \hat{v}_k(\hat{v}_{k-1}\cdots \hat{v}_2(v_1)) = \sum_{j=1}^n \hat{v}_k(\hat{v}_{k-1}\cdots \hat{v}_2(v_{1, j}))\frac{\partial}{\partial x_j}.$$ Using this notation, fundamental forms of $Z \subset \BP V$ can be described as follows.

\begin{itemize}
\item[(1)] Consider the map sending $v_{1,z} \otimes v_{2,z} \in T_z Z \otimes T_z Z$ to ${\rm d}_{z'} \xi (\hat{v}_1(\hat{v}_2))$ modulo $T_z Z$ where $\hat{v}_1$ (resp. $ \hat{v}_2$) is a local vector field on $\hat{Z}$ whose value at $z' \in \hat{z} \setminus 0$ is sent to $v_{1,z}$ (resp. $ v_{2,z}$) by ${\rm d}_{z'} \xi$. This determines a    homomorphism $${\rm FF}^2_{z,Z}: \Sym^2(T_z Z)
\longrightarrow T_z \BP V/ T_z Z$$ called the {\em second fundamental form} of $Z$ at $z$, independent of the choice of $z'$ and $\hat{v}_i$. There exists a Zariski open subset ${\rm Dom}({\rm FF}^3_Z) \subset Z$ such that the images of ${\rm FF}^2_{y,Z}, y \in {\rm Dom}({\rm FF}^3_Z),$ form a vector subbundle $$T^{(2)}Z \subset (T \BP V)|_{{\rm Dom}({\rm FF}^3_Z)}.$$  Thus we have a surjective vector bundle homomorphism $${\rm FF}_Z^2: \Sym^2(TZ) \to T^{(2)}Z/TZ \mbox{ on } {\rm Dom}({\rm FF}^3_Z) \subset Z.$$
\item[(2)] For $z \in {\rm Dom}({\rm FF}^{3}_Z),$ consider the map sending $v_{1,z} \otimes v_{2,z} \otimes  v_{3,z} \in T_z Z \otimes T_z Z \otimes T_z Z$  to ${\rm d}_{z'} \xi(\hat{v}_3(\hat{v}_{2} (\hat{v}_1)))$ modulo $T_z^{(2)}Z$, where $\hat{v}_{i}$ is a local vector field on $Z$ whose value at $z' \in \hat{z} \setminus 0$ is sent to  $v_{i,z}$ by ${\rm d}_{z'} \xi$ for each $i=1,2,3$. This determines a homomorphism
    $$ {\rm FF}^3_{z,Z}: \Sym^3(T_z Z)\longrightarrow T_z \BP V/T^{(2)}_z (Z),$$
called the {\em third fundamental form} of $Z$ at $z$, independent of the choices of $z'$ and $\hat{v}_i$. There exists a Zariski open subset ${\rm Dom}({\rm FF}^4_Z) \subset {\rm Dom}({\rm FF}^3_Z)$ such that the images of ${\rm FF}^3_{y,Z}, y \in {\rm Dom}({\rm FF}^4_Z)$ form a vector subbundle $$T^{(3)}Z \subset (T \BP V)|_{{\rm Dom}({\rm FF}^4_Z)}.$$ We obtain a surjective vector bundle homomorphism $${\rm FF}^3_Z : \Sym^3(TZ) \to T^{(3)}Z /T^{(2)} Z \mbox{ on } {\rm Dom}({\rm FF}^4_Z) \subset Z.$$
\item[(3)] Repeating the above construction inductively, we have  a homomorphism
    $$ {\rm FF}^k_{z,Z}: \Sym^k(T_z Z)\longrightarrow T_z \BP V/T^{(k-1)}_z (Z),$$ called the $k$-th {\em fundamental form} of $Z$ at $z$,
    where ${\rm Dom}({\rm FF}^{k}_Z) \subset {\rm Dom}({\rm FF}^{k-1}_Z)$ is a Zariski open subset
    such that the images of ${\rm FF}^{k-1}_{y,Z}, y \in {\rm Dom}({\rm FF}^k_Z),$
form a vector subbundle $$T^{(k-1)} Z \subset (T \BP V)|_{{\rm Dom}({\rm FF}^k_Z)}.$$ The quotient bundle $T^{(k)}Z / T^{(k-1)}Z$ on ${\rm Dom}({\rm FF}^{k+1}_Z)
\subset {\rm Dom}({\rm FF}^k_z)$ is denoted by $N^{(k)}_Z$ and called the $k$-th {\em normal space} of $Z$ with the fiber $N^{(k)}_{Z,z}$ called the $k$-th {\em normal space of  $Z$ at $z$}.  Thus we have a surjective vector bundle homomorphism $${\rm FF}_Z^k: \Sym^k(TZ) \to T^{(k)}Z / T^{(k-1)}Z = N^{(k)}_Z$$ on  ${\rm Dom}({\rm FF}^{k+1}_Z) \subset Z.$

\end{itemize}
We say that $Z \subset \BP V$ is {\em linearly nondegenerate} if it is not contained in any hyperplane in $\BP V$.
It is well-known that there exists $k \geq 1$ such that  $T^{(k)}_z Z = T_z \BP V$ for  any $z\in {\rm Dom}({\rm FF}^{k+1}_Z)$ if and only if $Z \subset \BP V$ is linearly nondegenerate.

\subsection{Jet spaces}\label{ss.j}

We recall the construction of jet spaces of functions for one independent and one dependent variables. In comparison with the references such as Section I.3 of \cite{BCG} or Chapter 4 of \cite{Ol}, we emphasize the projective bundle structure of the prolongation for the application in Section \ref{s.Goursat}.

Let $\Delta$ be a 1-dimensional complex manifold. Let $J^0 = \Delta \times \C$ be the trivial line bundle on $\Delta$ with the natural projection $p^0_{-1}: J^0 \to \Delta$ and set $\sJ^0 = TJ^0$ to be the tangent bundle.

Let $p^1_0: \BP \sJ^0 \to J^0$ be the projectivization of $\sJ^0$ and let  $\sJ^1 \subset T \BP \sJ^0$ be the prolongation of $\sJ^0$.
     Note that $\sJ^1$ is a contact structure on the 3-dimensional manifold  $\BP \sJ^0$ in the sense of Definition \ref{d.contact}.
      In fact, the contact form is $\vartheta^{J^0}$ of Example \ref{e.cotangent} under the natural isomorphism $\BP T J^0 = \BP T^* J^0$.
Let $p^2_1: \BP \sJ^1 \to \BP \sJ^0$ be the projectivization of $\sJ^1$ and let  $\sJ^2 \subset T \BP \sJ^1$ be the prolongation of $\sJ^1$.
Assuming we have defined a distribution $\sJ^k$ of rank 2 on a manifold $\BP \sJ^{k-1}$, we define  a distribution $\sJ^{k+1}$ of rank $2$ on $\BP \sJ^k$ as the prolongation of $\sJ^k$.

The {\em space $J^k$  of  $k$-jets} of functions on $\Delta$ is an open subset of $\BP \sJ^{k-1}$ defined as follows.
Assume that $\Delta$ has a coordinate function $t$ and let $u$ be a linear coordinate function on $\C$ such that $(t, u)$ is a coordinate system for $J^0$.
Let $J^1$ be the open subset $\BP \sJ^0 \setminus \BP T^{p^0_{-1}}.$ Regarding ${\rm d} t$ and ${\rm d} u$ as  holomorphic functions on $\sJ^0 = T J^0,$ we can view the meromorphic function $$u^{(1)} := \frac{{\rm d}u}{{\rm d} t}$$ on $\BP \sJ^0$ as a holomorphic function on $J^1$ defining  affine coordinates along the fibers of $J^1 \to J^0$.
Let $J^2$ be the open subset $\BP \sJ^1|_{J^1} \setminus \BP T^{p^1_0}.$
Then the meromorphic function $$u^{(2)} := \frac{{\rm d} u^{(1)}}{{\rm d} t}$$ on $\BP \sJ^1$ can be viewed as a holomorphic function on $J^2$ defining  affine coordinates along the fibers of $J^2 \to J^1$. Inductively define $J^{k+1}$ as the open subset $\BP \sJ^k|_{J^k} \setminus \BP T^{p^k_{k-1}}.$ Then the meromorphic function $$u^{(k+1)} := \frac{{\rm d} u^{(k)}}{{\rm d}t}$$ on $\BP \sJ^{k-1}$ can be viewed as a holomorphic function on $J^k$ defining affine coordinates along the fibers of $J^{k+1} \to J^k.$

A section of the fibration $J^{k+1} \to J^k$ over an open subset $U \subset J^k$ is given by $$u^{(k+1)} = F(t, u, u^{(1)}, \ldots, u^{(k)})$$
for a holomorphic function $F$ on $U$. This is equivalent to an ordinary differential equation of order $k+1$ in the independent variable $t$ and the dependent variable $u$ with initial conditions in $U$.

For $k > i \geq 1$, denote by $p^k_i: J^k \to J^i$ the composition $$p^{i+1}_i \circ p_{i+1}^{i+2} \circ \cdots \circ p^{k-1}_{k-2} \circ p^k_{k-1}.$$ It is easy to check from (\ref{e.pr}) (see also Theorem 6.5 of \cite{Mo}) that  the sequence $$\sJ^k \subset (p^k_{k-1})^{-1}\sJ^{k-1} \subset \cdots \subset (p^k_1)^{-1} \sJ^1 \subset T J^k$$ coincides with the sequence of  derived distributions on $J^k,$ both
$$\sJ^k \subset \partial^1 \sJ^k \subset \cdots \subset \partial^{k-1} \sJ^k \subset \partial^k \sJ^k$$ and
$$\sJ^k \subset \partial^{(1)} \sJ^k \subset \cdots \subset \partial^{(k-1)} \sJ^k \subset \partial^{(k)} \sJ^k.$$

\section{Distributions on the deformation spaces of unbendable rational curves}\label{s.vmrt}

 \begin{definition}\label{d.unbendable}
Let $C \subset M$ be a nonsingular projective rational curve on a complex manifold of dimension $n$. We say that $C$ is an {\em unbendable rational curve} if its normal bundle $N_C$ is isomorphic to \begin{eqnarray}\label{e.unbendable} N_C & \cong &  \sO_{\BP^1}(1)^{\oplus p} \oplus \sO_{\BP^1}^{\oplus (n-1-p)} \end{eqnarray} for some nonnegative integer $p$. The number $p+2$ is the {\em anti-canonical degree} of $C$, i.e., the degree of the line bundle $\wedge^n TM |_C$.
\end{definition}

By the deformation theory of rational curves, the family of nontrivial deformations of $C$ in $M$ has dimension $n-1+p = \dim H^0(C, N_C)$ and  all general deformations of $C$ in $M$ are also unbendable. It is convenient to formulate this in terms of the Douady space (or Hilbert scheme when $M$ is algebraic) of $M$ as follows. A good reference for Douady spaces is Section VIII.1 of \cite{GPR}.

\begin{definition}\label{d.sK} A {\em family of unbendable rational curves} on a complex manifold $M$ means  a connected open subset $\sK$ in the Douady space ${\rm Douady}(M)$ such that each member of $\sK$ is an unbendable rational curve. The
{\em anti-canonical degree} of $\sK$ is the anti-canonical degree of a member of $\sK$, namely, the integer $p+2$ in Definition \ref{d.unbendable}. We have $\dim \sK = n-1+p$.  Denote by
$$ \sK \stackrel{\rho}{\longleftarrow} \sU \stackrel{\mu}{\longrightarrow} M$$ the associated universal family morphisms. The morphism $\rho$ is a $\BP^1$-bundle and the morphism $\mu$ is a submersion
with the fiber dimension $p$. A point $y \in \sK$ corresponds to the unbendable rational curve $C_y := \mu(\rho^{-1}(y))$ in $M$. \end{definition}

 The following lemma is well-known from the deformation theory of rational curves (see page 58 of \cite{HM04}).

\begin{lemma}\label{l.tangents}
In Definition \ref{d.sK},  for a point $\alpha \in \sU$ with
$y = \rho(\alpha)$ and $x = \mu(\alpha)   \in C_y,$
we have  the following natural identifications.
\begin{itemize} \item[(i)] $T_y \sK = H^0(C_y, N_{C_y}).$
\item[(ii)] $T_{\alpha} \sU = H^0(C_y,  TM|_{C_y})/ H^0(C_y, TC_y \otimes {\bf m}_x).$
\item[(iii)] Under the identification in (ii), the subspace $T^{\rho}_{\alpha} \subset T_{\alpha} \sU$ tangent to the fiber of $\rho$ is identified with the subspace $$ H^0(C_y, T C_y)/H^0(C_y, T C_y \otimes {\bf m}_x) \subset H^0(C_y,  TM|_{C_y})/ H^0(C_y, TC_y \otimes {\bf m}_x)$$ induced by $TC_y \subset TM|_{C_y}$.
    \item[(iv)] Under the identification in (ii), the subspace $T^{\mu}_{\alpha}$ tangent to the fiber of $\mu$ is identified with $$  H^0(C_y, TM|_{C_y} \otimes {\bf m}_x)/ H^0(C_y, TC_y \otimes {\bf m}_x) = H^0(C_y, N_{C_y} \otimes {\bf m}_x).$$  \end{itemize}
       In particular, we have $T^{\mu}_{\alpha} \cap T^{\rho}_{\alpha} =0$ for any $\alpha \in \sU$ from $TM|_{C_y} \cong TC_y \oplus N_{C_y}$.  \end{lemma}

The space $\sK$ of unbendable rational curves carries a natural distribution:

\begin{definition}\label{d.sD}
For each member $C$ of $\sK$, let $N_C^+$ be the subbundle of the normal bundle $N_C$ corresponding to the $\sO_{\BP^1}(1)^{\oplus p}$ factor of (\ref{e.unbendable}), which is  independent of the choice of the isomorphism (\ref{e.unbendable}). Define a distribution $\sD \subset T \sK$ of rank $2p$ such that its fiber at $[C] \in \sK$ is
$$\sD_{[C]} := H^0(C, N_C^+) \subset H^0(C, N_C) = T_{[C]} \sK$$ via the natural identification in Lemma \ref{l.tangents} (i). \end{definition}

This distribution $\sD$ on $\sK$ is our main object of study. We relate it to some natural distributions on $\sU$ defined as follows.

\begin{definition}\label{d.sT} In Lemma \ref{l.tangents},
write $\sV := T^{\mu}$ and  $\sF := T^{\rho}$ on $\sU$. Set $\sT = \sV \oplus \sF$, which is a distribution on $\sU$ by Lemma \ref{l.tangents}.  Define meromorphic distributions
  $\sT^{k}, k \geq 0$ on $\sU$ inductively by  $\sT^0:= \sT$ and $$\sT^{k+1}:=[\sV, \sT^k]+\sT^k \mbox{ for each } k\geq 0.$$ Write $\sT_{k+1}:=\sT^{k+1}/\sT^k$ for each $k\geq 0$ and $\sT_0 := \sF $.
\end{definition}

The followings are straight-forward from the definition of $\sT^k$.

\begin{lemma}\label{l.psi}
In Definition \ref{d.sT}, there are  homomorphisms of sheaves $$ \psi_k: \sV\otimes\sT_{k} \rightarrow\sT_{k+1}, \ k\geq 0$$ on a Zariski open subset on $\sU$ such that $\psi_0(v\otimes f)=[v, f]$ modulo $\sT^0$  and $$\psi_k(v_{k+1}\otimes [v_{k}, [v_{k-1}, \cdots, [v_1, f]]])=[v_{k+1}, [v_{k},\cdots, [v_1, f]]]$$  for each $k\geq 1$,
where $f$ is a local section of $\sF$, and $v, v_1,\ldots, v_{k+1}$ are local sections of $\sV$. In the above equation, the repeated brackets representing elements of $\sT^k$ and $\sT^{k+1}$ are used to denote the corresponding elements of $\sT_k$ and $\sT_{k+1}$.
\end{lemma}

A key fact is the following result from Proposition 1 and Proposition 8 of \cite{HM04}.

\begin{proposition}\label{p.HM}
In Definitions \ref{d.sD} and \ref{d.sT}, the meromorphic distribution $\sT^1$ is a distribution on $\sU$ and satisfies $\sT^1 = \rho^{-1} \sD.$ \end{proposition}

\begin{proposition}\label{p.FJ} In the setting of Definition \ref{d.sT},
\begin{itemize}
\item[(i)] $[\sF, \sT^1]\subset\sT^1$; \item[(ii)]  $[\sF, \sT^2]\subset\sT^2$; and
\item[(iii)] If $[\sF, \sT^k]\subset\sT^k$ for $1 \leq k\leq 2m-1$ for some integer $m \geq 1$, then $[\sF, \sT^{2m}]\subset\sT^{2m}$. \end{itemize}
\end{proposition}

\begin{proof}
(i)  is a direct consequence of Proposition \ref{p.HM}. (ii) follows from (i) and (iii).

To prove (iii), we introduce the following notation. Fix  local sections $f$ of $\sF$ and $v_1,\ldots,v_m$ of $\sV$. Write  $$[v_1, \ldots, v_k, f] :=[v_1,[ v_2, [ \cdots [ v_k, f]\cdots]]]$$ which is a local section of  $\sT^k$, respectively.
The following lemma is immediate from Jacobi identity and $[\sV, \sV ] \subset \sV$.

\begin{lemma}\label{l.permutation}
For any $k \geq 1$ and any permutation $\sigma$ of the set $\{ 1, 2, \ldots, k\}$, $$[v_1, \ldots, v_k,f] - [v_{\sigma(1)}, \ldots, v_{\sigma(k)},f]$$ is a local section of $\sT^{k-1}$. \end{lemma}

Assuming that
 $[\sF, \sT^k]\subset\sT^k$ for all $1 \leq k\leq 2m-1,$  we have the following three lemmata.

\begin{lemma}\label{l.ij}
For any $0 \leq i \leq j\leq 2m-1$ satisfying $i+j\leq 2m-1$,
$$(\dagger)_{i,j} \ \ \ [[v_1,\ldots, v_i, f], [v_{i+1},\ldots, v_{i+ j}, f]] \mbox{ is a section of } \sT^{i+j}. $$ Here $[v_1,\ldots, v_i, f]$ means $f$ if $i=0$.
\end{lemma}

\begin{proof}
It is sufficient to prove
$$(\natural)_i: \mbox{ the statement } (\dagger)_{i, j} \mbox{ is true for all } i\leq j\leq 2m-1-i$$ for all $0\leq i < m$.

The statement $(\dagger)_{0,0}$ is immediate from $[\sF, \sF] \subset \sF$.
 When $i=0$ and $ 1 \leq j \leq 2m-1$, the statement $(\dagger)_{0,j}$ follows from the assumption $[\sF, \sT^j]\subset\sT^j$ for $j \leq 2m-1$. This proves  $(\natural)_0$.

Now fix $1 \leq i < m$ and assume that $(\natural)_k$ is true for all $0 \leq k < i$. Take any $i\leq j\leq 2m-1-i$.
Then
$$
[[v_1, \ldots, v_i, f], [v_{i+1}, \ldots, v_{i+j}, f]]$$ is equal to
$$ [[v_1, v_{i+1}, \ldots, v_{i+j},f], [v_2, \ldots, v_i, f]] + [v_1, [[v_2, \ldots, v_i,f ], [v_{i+1}, \ldots, v_{i+j}, f]]].$$ The first term belongs to  $\sT^{i+j}$ by the induction hypothesis $(\dagger)_{i-1, j+1}$  and the second term belongs to $[\sV, \sT^{i+j-1}] \subset \sT^{i+j}$ by the induction hypothesis $(\dagger)_{i-1, j}.$ This proves $(\natural)_i$ and the lemma. \end{proof}

\begin{lemma}\label{l.split}
For any $1 \leq k \leq 2m,$
$$[[v_1, \ldots, v_{k-1},f], [v_k, \ldots, v_{2m}, f]] \equiv - [[v_1, \ldots, v_k,f], [v_{k+1}, \ldots, v_{2m}, f]] $$ modulo $\sT^{2m}$.
Here $[v_1, \ldots, v_{k-1},f]$ means $f$ if $k=1$ and $[v_{k+1}, \ldots, v_{2m}, f]$  means $f$ if $k=2m$.   \end{lemma}

\begin{proof}
By Jacobi identity,
\begin{eqnarray*} \lefteqn{ [[v_1,\ldots, v_{k-1},f], [v_k, \ldots, v_{2m},f]] } \\ & = &
[[v_1, \ldots, v_{k-1},f], [v_k, [v_{k+1}, \ldots, v_{2m}, f]]] \\ &=&
 -[[v_k, [v_1, \ldots, v_{k-1},f]], [v_{k+1}, \ldots, v_{2m}, f]]  \\ && + [v_k, [[v_1, \ldots, v_{k-1},f], [v_{k+1}, \ldots, v_{2m}, f]]].\end{eqnarray*}
The last term  is a local section of $\sT^{2m}$ by Lemma \ref{l.ij}. The second last term  is
$$ - [[v_k, v_1, \ldots, v_{k-1}, f], [v_{k+1}, \ldots, v_{2m}, f]] $$
which is equal to $-[[v_1, \ldots, v_k, f],[ v_{k+1}, \ldots, v_{2m}, f]]$ modulo $\sT^{2m}$ by Lemma \ref{l.permutation}.
\end{proof}

\begin{lemma}\label{l.3}
For any $1 \leq k \leq 2m$, we have  $$[f, [v_1,\ldots, v_{2m},f]] \equiv (-1)^k [[v_1,\ldots, v_k, f], [v_{k+1}, \ldots, v_{2m}, f]]$$ modulo $\sT^{2m}$. Here $[v_{k+1}, \ldots, v_{2m}, f]$ means $f$ if $k=2m$. \end{lemma}

\begin{proof}
This follows from Lemma \ref{l.split} by induction on $k$. \end{proof}

  Lemma \ref{l.3} with $k=2m$ says
  $$[f, [v_1, \ldots, v_{2m}, f]] \equiv (-1)^{2m} [[v_1, \ldots, v_{2m}, f], f]$$ modulo $\sT^{2m}$. Thus $[f, [v_1, \ldots, v_{2m}, f]] $ is a local section of $\sT^{2m}.$   Since elements of the form $[v_1, \ldots, v_{2m},f]$ span $\sT^{2m}$ locally, we obtain Proposition \ref{p.FJ} (iii).  \end{proof}

\begin{remark}
 Proposition \ref{p.FJ} (i) and (ii) cannot be extended further. In fact,  there are cases where $[\sF, \sT^3] \not\subset \sT^3$, as explained in Sections \ref{s.Cartan} and \ref{s.quartic} below.
\end{remark}

\begin{corollary}\label{c.T2}
In setting of Definition 3.5, we have $$\partial^{(2)}\sT^0=\partial^2\sT^0=\sT^2=\rho^{-1}(\partial\sD).$$
\end{corollary}

\begin{proof}
By definition, we have $$\sT^2=[\sV, \sT^1]+\mathcal{T}^1\ \subset \ \partial^{(2)}\sT^0 \ \subset \ \partial^2\sT^0 \ = \ \partial\sT^1.$$ For  a local section $f$ of $\mathcal{F}$ and local sections $v_1, v_2$ of $\mathcal{V}$,  we have
\begin{eqnarray*}
[[v_1, f], [v_2, f]]=[[v_1, [v_2, f]], f]+[v_1, [f, [v_2, f]]].
\end{eqnarray*}
As $[[v_1, [v_2, f]], f]$ and $[v_1, [f, [v_2, f]]]$ are both local sections of $\mathcal{T}^2$ by Proposition \ref{p.FJ} (i) and (ii),  so is $[[v_1, f], [v_2, f]]$. Since elements of the form $[v_i, f]$ span $\mathcal{T}^1$ locally, we have $\partial\mathcal{T}^1\subset\mathcal{T}^2$ and consequently
$$\mathcal{T}^2=\partial^{(2)}\mathcal{T}^0=\partial^2\mathcal{T}^0.$$
Using $\mathcal{T}^1=\rho^{-1}\mathcal{D}$ from Proposition \ref{p.HM}, we obtain $\partial^2\mathcal{T}^0=\partial\mathcal{T}^1=\rho^{-1}(\partial\mathcal{D})$ by (\ref{e.dist}).
\end{proof}

We recall the definition of varieties of minimal rational tangents. Originally (as in \cite{HM}, \cite{Hw01}), they are defined for minimal rational curves on uniruled projective manifolds. But essentially the same definition works for any unbendable rational curve in complex manifolds.

\begin{definition}\label{d.vmrt}
In the terminology of Lemma \ref{l.tangents}, define the {\em tangent morphism } $\tau: \sU \to \BP T M$ by setting $$\tau(\alpha) := [T_x C_y] \in \BP T_x M.$$ Both $\tau$ and the restriction $$\tau_x: \mu^{-1} (x) \to \BP T_x M$$ for each $x \in M$ are immersions (e.g. by Proposition 1.4 of \cite{Hw01}).
The image $\sC_x$ of $\tau_x$ is called the {\em variety of minimal rational tangents} (VMRT in abbreviation) of $\sK$ at $x \in M$.
For $\alpha \in \mu^{-1}(x)$, the $\tau_x$-image of the germ of $\mu^{-1}(x)$ at $\alpha$ is denoted by $\sC_{\alpha} \subset \BP T_x M$. Since $\tau_x$ is immersive, the germ $\sC_{\alpha}$ is a submanifold and an irreducible component of the germ of $\sC_x$ at $\tau_x(\alpha)$.
\end{definition}

The next proposition is a generalization of Proposition 2 of \cite{HM04}.

\begin{proposition}\label{p.sT}
In the setting of Lemma \ref{l.psi}, the following holds.   \begin{itemize}

%\item[(i)] The line bundle $T^{\rho}$ on $\mu^{-1}(x)$ is isomorphic to $\tau_x^* \sO(-1)_{\BP T_x M}$.

\item[(i)] The homomorphism $\psi_0$  is an isomorphism between $\sV \otimes \sF$ and $\sT_1$.

\item[(ii)] Pick $\alpha \in \sU$ such that
$z:= \tau(\alpha)$ is in ${\rm Dom}({\rm FF}^k_{\sC_{\alpha}})$ for some $k\geq 2$.
Then for a local section $f$ of $\sF$ and local sections $v_1, \ldots, v_k$ of $\sV$ near $\alpha$, the value of the vector field $[v_k, [v_{k-1},\cdots, [v_1, f]]]$ at $\alpha$
satisfies $${\rm d}\mu ( [v_k, [v_{k-1},\cdots, [v_1, f]]]_{\alpha})
 \equiv {\rm d}\mu (f_{\alpha}) \otimes {\rm FF}^k_{z, \sC_{\alpha}}(v_1,\ldots, v_k)$$ modulo $T^{(k-1)}_z \sC_{\alpha}$.

\item[(iii)] In the setting of (ii), we have a natural identification $$\sT_{k, \alpha} =\hat{z}\otimes N^{(k)}_{\sC_{\alpha},z}$$ where  $N^{(k)}_{\sC_{\alpha},z}$ is the $k$-th normal space of $\sC_{\alpha} \subset\BP T_x M$ at the point $z\in\sC_{\alpha}$. In particular, the homomorphism $\psi_j$ of Lemma \ref{l.psi} is surjective on ${\rm Dom}({\rm FF}^k_{\sC_{\alpha}})$ for all $ 1\leq j \leq k$.

\item[(iv)] If $\sC_{\alpha}$ is linearly nondegenerate in $\mathbb{P}T_x M, x = \mu(\alpha)$ for a general $\alpha \in \sU$, then $T\sU=\sT^k$ on a Zariski open subset of $\sU$ for a sufficiently large integer $k$.

\item[(v)] In the setting of (iv), the distribution $\mathcal{D}$ on $\mathcal{K}$ is bracket-generating.
    \end{itemize}
\end{proposition}

\begin{proof}
(i) is from Section 2 of \cite{HM04}.

To prove (ii), we modify the computation in the proof of Proposition 2 in \cite{HM04} as follows.
Fix a local coordinate system  $x_1,\cdots, x_n$ in a neighborhood of  $x \in M$. Then $\lambda_1:=dx_1,\ldots, \lambda_n:=dx_n$ give linear coordinates in the vertical directions of $TM \rightarrow M$.

Via the immersion $\tau$, let us identify a neighborhood $U$ of $\alpha$ in $\sU$ with a submanifold of $\BP TM$  such that  $\alpha$ is identified with $z$ and $\mu$ agrees with the restriction of the natural projection $\BP TM \to M$.
Let $\xi: TM\setminus(0\mbox{-section})\rightarrow \BP TM$ be the natural $\C^*$-bundle. We can choose a point $z' \in \xi^{-1}(z)$, a neighborhood $O$ of $z'$ in $TM$ and local holomorphic functions $f_i, v_{i,j}$ on $O$ with $1\leq i \leq n$ and $ 1 \leq j \leq k$   such that the  vector fields on $O$ $$\hat{f}=\sum_{j=1}^n f_j\frac{\partial}{\partial\lambda_j}+\sum_{j=1}^n\lambda_j\frac{\partial}{\partial x_j}, \ \hat{v}_i=\sum_{j=1}^n v_{i, j}\frac{\partial}{\partial\lambda_j}$$ are tangent to $\xi^{-1}(U)$   and their images under  ${\rm d} \xi $ induce  $f$ and $v_i$ on $U$.
  Direct calculations show that $$[\hat{v}_1, \hat{f}]=\sum_{j=1}^n v_{1, j}\frac{\partial}{\partial x_j}  $$ modulo $\frac{\partial}{\partial \lambda_1}, \ldots, \frac{\partial}{\partial \lambda_n}$ and $$[\hat{v}_m, [\hat{v}_{m-1},\cdots, [\hat{v}_1, \hat{f}]]]=\sum_{j=1}^n \hat{v}_m(\hat{v}_{m-1}\cdots \hat{v}_2(v_{1, j}))\frac{\partial}{\partial x_j} $$  modulo $\frac{\partial}{\partial \lambda_1}, \ldots, \frac{\partial}{\partial \lambda_n}$ for each $2\leq m\leq k$. This proves (ii).

(iii) and (iv) follow from (ii). Since $\mathcal{T}^1=\rho^{-1}\mathcal{D}$ by Proposition \ref{p.HM}, we have $\mathcal{T}^{k+1}\subset\partial^k\mathcal{T}^1=\rho^{-1}\partial^k\mathcal{D}$ for each $k\geq 1$. Thus (v) follows from (iv).
\end{proof}

The next proposition is well-known. (i) is from Proposition 2.3 of \cite{Hw01} and the rest follows easily from (i) (see Proposition 3.1 of \cite{HH} or the proof of Propositions 2.2 and 3.1 of \cite{Mk}).

\begin{proposition}\label{p.FFvmrt}
In the setting of  Lemma \ref{l.psi} and Definition \ref{d.vmrt},
for a point $y \in \sK$, define $$\sC_y:= \cup_{\alpha \in \rho^{-1}(y)} \sC_{\alpha} \ \subset \BP TM|_{C_y}.$$ It is a complex manifold with a submersion $\pi_y: \sC_y \to C_y$.
 \begin{itemize}
 \item[(i)] The relative tangent bundle of $\pi_y: \sC_y \to C_y$ restricted to  $\rho^{-1}(y) \cong \BP^1$ is isomorphic to $\sO_{\BP^1}(-1)^{\oplus p}.$
 \item[(ii)] The relative normal bundle of $\sC_y \subset \BP TM|_{C_y}$ restricted to $\rho^{-1}(y)$ is isomorphic to $\sO_{\BP^1}(-2)^{\oplus (n-1-p)}$ and the relative second fundamental forms of $\sC_y \subset \BP TM|_{C_y}$ along $\rho^{-1}(y)$ are given by a section of a vector bundle on $\rho^{-1}(y)$ isomorphic to $$\Hom (\Sym^2\sO_{\BP^1}(-1)^{\oplus p}, \sO_{\BP^1}(-2)^{\oplus (n-1-p)}) = \sO_{\BP^1}^{\oplus (p^2+p)(n-1-p)/2}.$$
 \item[(iii)] If a point $\alpha \in \rho^{-1}(y)$ is in ${\rm Dom}({\rm FF}^3_{\sC_{\alpha}})$, then any point $\beta \in \rho^{-1}(y)$ is in ${\rm Dom}({\rm FF}^3_{\sC_{\beta}}).$ In this case, the relative second normal spaces of $\sC_y \subset \BP TM|_{C_y}$ along $\rho^{-1}(y)$ determine a vector bundle on $\rho^{-1}(y)$ isomorphic to $\sO_{\BP^1}(-2)^{\oplus q}$ for some nonnegative integer $q\leq n-1-p.$
 \item[(iv)] When a point $\alpha \in \rho^{-1}(y)$ is in ${\rm Dom}({\rm FF}^3_{\sC_{\alpha}}),$ the relative third fundamental forms of $\sC_y \subset \BP TM|_{C_y}$ along $\rho^{-1}(y)$ are given by a section of a vector bundle on $\rho^{-1}(y)$ isomorphic to $$\Hom (\Sym^3\sO_{\BP^1}(-1)^{\oplus p}, \sO_{\BP^1}(-2)^{\oplus (n-1-p-q)}) = \sO_{\BP^1}(1)^{\oplus r}$$ for $r = \frac{1}{6}(p^3+3p^2+2p)(n-1-p-q)$. \end{itemize} \end{proposition}

The above results say that the properties of $\sD$ are intimately related to the projective geometry of VMRT, especially their fundamental forms. Many different types of fundamental forms of VMRT can arise depending on the geometry of the unbendable rational curves and it is difficult to obtain uniform results. We have to restrict to specific classes of unbendable rational curves to proceed further.

\begin{definition}\label{d.bg}
A family $\sK$ of unbendable rational curves on a complex manifold is said to be {\em bracket-generating} if the distribution $\sT$  is a bracket-generating distribution on $\sU$. By Proposition \ref{p.HM}, this is equivalent to saying that $\sD$ in Definition \ref{d.sD} is a bracket-generating distribution on $\sK$. \end{definition}

The next proposition shows that to study germs of unbendable rational curves, we may concentrate on bracket-generating families.

\begin{proposition}\label{p.bg}
Suppose $\sK$ is a family of unbendable rational curves on a complex manifold $M$, which is not bracket-generating. Then a general member
$C\subset M$ of $\sK$ has a neighborhood $C \subset U \subset M$ equipped with a submersion $\eta: U \to B$ to a complex manifold $B$ such that $\eta(C)$ is a point $b \in B$ and the family $\sK_b,$ which consists of members of $\sK$ representing rational curves contained  in $\eta^{-1}(b),$ is bracket-generating.  In particular, the curve $C$ is an unbendable rational curve in the lower-dimensional complex manifold $\eta^{-1}(b)$. \end{proposition}

\begin{proof}
To prove Proposition \ref{p.bg}, we may replace $\sK$ by a neighborhood of a general member of $\sK$. In particular, we may assume that $\mu$ has connected fibers and $\sC_x$ is irreducible for  a general point $x \in M$. We define a meromorphic distribution $\sE \subset TM$ such that its fiber $\sE_x \subset T_xM$ at a general point $x\in M$ is the linear span of  $\sC_x$ in $T_x M$.
Then $\sT^k = \mu^{-1} \sE$ on an open subset of $\sU$ for sufficiently
large $k$. It follows that $\sE$ is not bracket-generating on $M$ because $\sT$ is not bracket-generating on $\sU$.
Thus $\partial^{m} \sE$ becomes involutive for some $m \geq 1$ and has rank strictly smaller than $\dim M.$ We have  a closed analytic subset $Z \subset M$ of codimension $\geq 2$ such that $\partial^{m} \sE$ is an involutive distribution on $M \setminus Z$.
Then a general member $C$ of $\sK$ is disjoint from $Z$ by Lemma 2.1 of \cite{Hw01}. Then the leaves of the foliation defined by $\partial^{m} \sE$ induces a submersion $\eta:U \to B$ in a neighborhood $U$ of $C$ such that $C$ is contained in a fiber of $\eta$. It is clear that $C$ is an unbendable rational curve on the fiber and $\sD|_{\sK_b} \subset T \sK_b$ is bracket-generating.  \end{proof}

For the rest of the paper, we concentrate on unbendable rational curves with $p=1,$ which is exactly when $\sD$ is a distribution of rank $2$. In this case, we have the following general result.

\begin{proposition}\label{p.rank2}
For a family of unbendable rational curves of anti-canonical degree $3$ on $M$, let $\pi: \BP \sD \to \sK$ be the projectivization of the vector bundle $\sD$ of rank 2 on $\sK$. \begin{itemize} \item[(i)] Let $\iota: \sU \to \BP \sD$ be the map sending   $\alpha \in \sU$ with $y= \rho(\alpha)$ and  $x = \mu(\alpha)$ to the 1-dimensional subspace of $H^0(C_y, N_{C_y}^+) = \sD_y$ given by $H^0(C_y, N_{C_y}^+ \otimes {\bf m}_x)$. Then $\iota$ is a biholomorphic map satisfying $\rho= \pi \circ \iota.$ \item[(ii)]  The image ${\rm d} \iota(T^{\mu})$ is a line subbundle of ${\rm pr}(\sD) \subset \pi^{-1} \sD$ satisfying ${\rm d} \iota(T^{\mu}) \cap T^{\pi} =0$, i.e., it splits the exact sequence $$ 0 \to T^{\pi} \to {\rm pr}(\sD) \to \sO_{\BP \sD} (-1) \to 0.$$ \end{itemize}
Conversely, given a distribution $D$ of rank $2$ on a complex manifold $Y$, let $\pi: \BP D \to Y$ be its projectivization and let $\sV \subset {\rm pr}(D) $ be a line subbundle on a neighborhood $U$ of a fiber of $\pi$ that splits $$ 0 \to T^{\pi} \to {\rm pr}(D) \to \sO_{\BP D} (-1) \to 0$$ on $U$. Then after shrinking $U$ if necessary, we obtain   a submersion $\mu: U \to M$ to a complex manifold $M$ whose fibers are the leaves of $\sV$ such that the images of the fibers of $\pi$ under $\mu$ define a family of unbendable rational curves of anti-canonical degree 3 on $M$ and $D$ can be identified with the distribution $\sD$  associated with the family. \end{proposition}

\begin{proof} Let $n$ be the dimension of $M$.
(i) follows easily from $N_{C_y} \cong \sO_{\BP^1}(1) \oplus \sO_{\BP^1}^{\oplus (n-2)}.$ (ii) follows from Lemma \ref{l.tangents} and Proposition \ref{p.FFvmrt} (i).
To see the converse statement, the curve $C_y = \mu(\pi^{-1}(y))$ corresponding to a point $y \in Y$ has the normal bundle $N_{C_y}$ of degree 1 which is generated by global sections. Thus it is isomorphic to $\sO_{\BP^1}(1) \oplus \sO_{\BP^1}^{\oplus (n-2)}$, implying that $C_y$ is unbendable. This proves that  a neighborhood of $y$ in $Y$ can be regarded as the family of unbendable rational curves in an $n$-dimensional complex manifold $M$. \end{proof}

\section{Unbendable rational curves of Goursat type}\label{s.Goursat}

\begin{definition}\label{d.Goursat}
A  distribution $D$ of rank 2 on an $n$-dimensional complex manifold $Y$ is {\em a Goursat distribution} if $\partial^{k} D$ is a distribution of rank $k+2$  for each $1 \leq k \leq n-2.$
A family $\sK$ of unbendable rational curves of anticanonical degree 3 in an $n$-dimensional complex manifold  is said to be {\em of Goursat type}, if the distribution $\sD$ on  $\sK$ in Definition \ref{d.sD} is a Goursat distribution on its regular locus. By Proposition \ref{p.HM}, this is equivalent to saying that the distribution $\sT$ on $\sU$ is a Goursat distribution on its regular locus. \end{definition}

An example of Goursat distribution is the distribution $\sJ^n$ on the jet space $J^n$ of Subsection \ref{ss.j}. The  following result of \'{E}. Cartan (Theorem 6.5 of \cite{Mo}) says that any Goursat distribution is isomorphic to $\sJ^n$ at general points.

\begin{theorem}\label{t.Cartan}
Let $D \subset TY$ be a Goursat distribution on a manifold $Y$ of dimension $n \geq 3$. Then at a general point $y\in Y$, there exist a neighborhood $y \in U \subset Y$ and an open subset $O \subset J^{n-2}$ in the jet space of dimension $n$ (see Subsection \ref{ss.j}) with a biholomorphic map $\varphi: U \to O$ such that ${\rm d} \varphi (D|_U) = \sJ^{n-2}|_O$. \end{theorem}

A direct consequence of Proposition \ref{p.rank2} and Theorem \ref{t.Cartan} is the following.

\begin{theorem}\label{t.Goursat}
Let $\sK$ be a family of unbendable rational curves of Goursat type on a complex manifold $M$ of dimension $n$ with the universal morphisms $\sK \stackrel{\rho}{\leftarrow} \sU \stackrel{\mu}{\rightarrow} M$.
By Theorem \ref{t.Cartan},  a general point of $\sK$ has a neighborhood $U$ biholomorphic to an open subset $O \subset J^{n-2}$ such that $\sD \subset T\sK$ corresponds to $\sJ^{n-2}$ on $O$. Then there is a biholomorphism $\phi$ with a commutative diagram $$ \begin{array}{ccc} \rho^{-1}(U) & \stackrel{\phi}{\longrightarrow} & \BP \sJ^{n-2}|_O \\ \downarrow & &  \downarrow \\ U & \stackrel{\varphi}{\longrightarrow} & O \end{array} $$ such that ${\rm d} \phi$ sends  $T^{\mu}|_{\rho^{-1}(U)}$  to a line subbundle of $\sJ^{n-1}$ on $\BP \sJ^{n-2}|_O$ that splits  the  exact sequence $$0 \to T^{p^{n-1}_{n-2}} \to \sJ^{n-1} \to \sO_{\BP \sJ^{n-2}}(-1) \to 0$$ over the open subset $(p^{n-1}_{n-2})^{-1}(O) \subset \BP \sJ^{n-2}$, where $p^{n-1}_{n-2}: \BP \sJ^{n-2}|_{J^{n-2}} \to J^{n-2}$ is the natural projection.  Conversely, given a line subbundle $\sV$ of $\sJ^{n-1}$ over an  open subset $O \subset J^{n-2}$ that splits the above exact sequence, let $M$ be the space of leaves of $\sV$ in a neighborhood of a fiber of $p^{n-1}_{n-2}.$ Then the images of the fibers of $p^{n-1}_{n-2}$ in $M$ give a family of unbendable rational curves of Goursat type on $M$. \end{theorem}

The following proposition gives an interesting class of examples of unbendable rational curves of Goursat type.

\begin{proposition}\label{p.flat}
Let $Z \subset \BP^{n-1}$ be a submanifold of dimension 1. Regard $\BP^{n-1} $ as a hyperplane in $\BP^n$ and let $\beta: M \to \BP^n$ be the blowup of $\BP^n$ along $Z$. Let $\sK$ be the set of  lines on $\BP^{n}$ which intersect $Z$ and are not contained in $\BP^{n-1}.$ Regard $\sK$ as the parameter space of the strict transforms to $M$ of such lines. Then $\sK$ is a family of  unbendable rational curves of anti-canonical degree $3$ on $M$. If furthermore $Z$ is linearly nondegenerate in $\BP^{n-1}$, then $\sK$ is of Goursat type. \end{proposition}

\begin{proof}
A line in $\BP^n$ belonging to $\sK$ intersects $Z$ transversally and its normal bundle  in $\BP^n$ is isomorphic to $\sO(1)^{\oplus (n-1)}$. Thus its strict transform in $M$ is a rational curve with normal bundle $\sO(1) \oplus \sO^{\oplus (n-2)}$, namely, an unbendable rational curve of anti-canonical degree $3$.

Now assume that $Z$ is linearly nondegenerate in $\BP^{n-1}$ and let us show that $\sK$ is of Goursat type.
The linear nondegeneracy of $Z \subset \BP^{n-1}$ implies \begin{eqnarray}\label{e.rank} {\rm rank} (\sT^i) &=& i+2 \mbox{ for } 1 \leq i \leq n-1\end{eqnarray} by Proposition \ref{p.sT}.

Let $M_o \subset M$ be the Zariski open subset given by $\beta^{-1}(\BP^n \setminus \BP^{n-1})$ and fix  affine coordinates $(x_1, \ldots, x_n)$ on $M_o$ induced by an inhomogeneous coordinate system on $\BP^n \setminus \BP^{n-1}$. Let $\lambda_1 = {\rm d} x_1, \ldots, \lambda_n= {\rm d} x_n$ be the fiber coordinates on $TM_o$.
The members of $\sK$ intersect $M_o$ on affine lines in the direction of $Z \subset \BP^{n-1}$ via the induced trivialization $\BP TM_o \cong \BP^{n-1} \times M_o$.
In other words, the submanifold $\sC:= \tau(\mu^{-1}(M_o)) \subset \BP TM_o$ is isomorphic to the product $Z \times M_o \subset \BP^{n-1} \times M_o = \BP TM_o$.
In terms of the coordinates $$( \lambda_1, \ldots, \lambda_n, x_1, \ldots, x_n)$$ on $TM_o$,  the vector field $$\hat{f}:= \lambda_1 \frac{\partial}{\partial x_1} + \cdots + \lambda_n \frac{\partial}{\partial x_n}$$ on $TM_o$ restricted to the cone $\hat{\sC} \subset TM_o \setminus \mbox{(0-section)}$ gives a section of $\xi^{-1} \sF$ where $\xi: \hat{\sC}  \to \sC$ is the natural $\C^*$-bundle. Local vector fields
on $Z$ can be extended to local vector fields on $\sC \cong Z \times M_o$. Thus fibers of $\xi^{-1} \sV$ on $\hat{\sC}$ are spanned by  restrictions to $\hat{\sC}$  of local vector fields of the form
$$a_1 \frac{\partial}{\partial \lambda_1} + \cdots + a_n \frac{\partial}{\partial \lambda_n}$$
where $a_1, \ldots, a_n$ are suitable holomorphic functions in $\lambda_1, \ldots, \lambda_n$ defined on some open subset in $\C^{n}$ such that they represent local vector fields tangent to  $\hat{Z} \subset \C^n$.
It follows that $\sT^1 = \sV + \sF + [\sV, \sF]$ is  generated by $\sV, \sF$ and local vector fields
of the form $$[ a_1 \frac{\partial}{\partial \lambda_1} + \cdots + a_n \frac{\partial}{\partial \lambda_n}, \hat{f}] =  a_1 \frac{\partial}{\partial x_1} + \cdots + a_n \frac{\partial}{\partial x_n}$$
where $a_1, \ldots, a_n$ are suitable local holomorphic functions in $\lambda_1, \ldots, \lambda_n$.
By induction, we see that $\sT^k$ for any $k$ is generated by $\sV, \sF$ and local vector fields
of the form $$b_1 \frac{\partial}{\partial x_1} + \cdots + b_n \frac{\partial}{\partial x_n}$$
where $b_1, \ldots, b_n$ are suitable local holomorphic functions in $\lambda_1, \ldots, \lambda_n$ and $\sum_{i=1}^n b_i \frac{\partial}{\partial \lambda_i}$ is a local section of $T^{(k)} \hat{Z}.$
This implies that  $$[\sF, \sT^i] \subset \sT^i \mbox{ and } [\sT^i, \sT^j] \subset \sT^{j+1} \subset \sT^{i+j} \mbox{ for any } 1 \leq i \leq j.$$ From
$\sT^1 = \rho^{-1}\sD$, we have $$\rho^{-1}(\partial^{(i)} \sD )= \sT^{i+1} = \rho^{-1}( \partial^{i} \sD).$$ Combined with $(\ref{e.rank})$, we see that the growth vector of $\sD$ is $(2, 3, \ldots, n)$ and $\sD$ is a Goursat distribution on its  regular locus. \end{proof}

\begin{definition}\label{d.ode}
Let $O \subset J^{n-2}$ be a connected open subset and let  $$\sV \subset \sJ^{n-1}|_{(p^{n-1}_{n-2})^{-1}(O)}$$ be a line subbundle  splitting the exact sequence in  Theorem \ref{t.Goursat} on $(p^{n-1}_{n-2})^{-1}(O) \subset \BP \sJ^{n-2}$.
The line bundle $\sV$ determines a section of $p^n_{n-1}: \BP \sJ^{n-1} \to \BP \sJ^{n-2}$ on $(p^{n-1}_{n-2})^{-1}(O)$. This section restricted to  the open subset $J^{n-1} \cap (p^{n-1}_{n-2})^{-1}(O)$ determines an ordinary differential equation
\begin{eqnarray}\label{e.ode} u^{(n)} &=& F(t, u, u^{(1)}, \ldots, u^{(n-1)}) \end{eqnarray}
for some holomorphic function $F$ on $J^{n-1} \cap (p^{n-1}_{n-2})^{-1}(O) $
as explained in Subsection \ref{ss.j}.
If $O$ and $\sV$ arise from a neighborhood of a general point $y \in \sK$ of a family of unbendable rational curves of Goursat type as described in Theorem \ref{t.Goursat}, we say that (\ref{e.ode}) is an {\em ODE associated with} $\sK$. It is not unique, depending on the choice of the biholomorphism $\varphi$. \end{definition}

  What kind of ordinary differential equations are associated with unbendable rational curves in Definition \ref{d.ode}? To answer it, we look at one special example first.

\begin{proposition}\label{p.special}
Let $Z \subset \BP^{n-1}$ be the rational normal curve of degree $n-1 \geq 2$, in other words, the set of pure symmetric tensors in $\BP \Sym^{n-1}(\C^2)$, and let $\sK$ be the family of unbendable rational curve of Goursat type determined by $Z \subset \BP^{n-1}$ in Proposition \ref{p.flat}.
Then we can choose $\varphi$ in Theorem \ref{t.Goursat} such that  $u^{(n)} =0$ is an ODE associated with $\sK$. \end{proposition}

\begin{proof}
The semi-direct product $(\Sym^{n-1} \C^2) \times {\rm GL}(2) =   \C^n \times {\rm GL}(2)$ acts on the configuration $$Z \subset \BP^{n-1} \subset \BP^n=\mathbb{P}(\mathbb{C}\oplus\Sym^{n-1}\mathbb{C}^2)$$ such that  the  kernel of the action is the finite subgroup
$$\{(0, a)\in\mathbb{C}^n\times{\rm GL}(2)\mid a\in\mathbb{C}^*{\rm Id}_{\C^2} \mbox{ and } a^{n-1}={\rm Id}_{\C^2} \}.$$
This action can be lifted to the blowup $M \to \BP^n$ along $Z$ with an induced action on $\sK$. It follows that $\dim \mathfrak{a}\mathfrak{u}\mathfrak{t}(\sD) \geq n+4$.
It is a classical result (for example, see Theorem 6.44 of \cite{Ol}) that any ordinary differential equation of order $n \geq 4$ with $(n+4)$-dimensional Lie algebra of infinitesimal automorphisms is equivalent (up to a choice of $\varphi$) to $u^{(n)} =0.$ This completes the proof when $n\geq 4$.

When $n =3$,  the curve $Z\subset \BP^2$ is a conic and the strict transform of $\BP^2 \subset \BP^3$ under the blowup $M \to \BP^3$ can be contracted to a point yielding a birational morphism from   $M$ to the nonsingular quadric hypersurface $\Q^3 \subset \BP^4$ (e.g. see Lecture 22 of \cite{Ha}, page 288).   Then $\sK$ can be viewed as a family of unbendable rational curves on $\Q^3$ and the infinitesimal automorphisms of $\Q^3$ induce infinitesimal automorphisms of $\sD$. It follows that $$\dim \mathfrak{a}\mathfrak{u}\mathfrak{t}(\sD) \geq \dim \mathfrak{s}\mathfrak{o}(\C^5) = 10.$$
 We have the classical result (see Theorem 6.44 of \cite{Ol} again)  that any third order ordinary differential equation with  $10$-dimensional Lie algebra of infinitesimal (contact) automorphisms is equivalent (up to a choice of $\varphi$) to $u^{(3)}=0.$  This completes the proof when $n=3$. \end{proof}

We have the following general result. This can be viewed as a generalization of   Theorem 3.1 in \cite{Hi} which corresponds to the case $n=2$.

\begin{theorem}\label{t.ode}
Let $\sK$ be a family of unbendable rational curves of Goursat type in a complex manifold of dimension $n \geq 3$. Then there exists an ODE associated with $\sK$ of the form \begin{eqnarray}\label{e.ode} u^{(n)} &=& a_3 (u^{(n-1)})^3 + a_2 (u^{(n-1)})^2 + a_1 u^{(n-1)} + a_0\end{eqnarray} where $a_0, \ldots, a_3$ are holomorphic functions of $t, u, u^{(1)}, \ldots, u^{(n-2)}$ on a domain in $J^{n-2}$.
Conversely, any ODE of this type gives rise to a family of unbendable rational curves of Goursat type in some complex manifold of dimension $n$. \end{theorem}

\begin{proof}
By Theorem \ref{t.Goursat}, a family of unbendable rational curves of Goursat type gives a splitting of the exact sequence $$0 \to T^{ p^{n-1}_{n-2}} \to \sJ^{n-1} \to \sO_{\BP \sJ^{n-2}}(-1) \to 0$$ over an open subset $O$ in $J^{n-2}$. The unbendable curves in Proposition \ref{p.special} give another splitting of this exact sequence. Their difference is a section of the line bundle $T^{ p^{n-1}_{n-2}} \otimes \sO_{\BP \sJ^{n-2}}(1) $ on $(p^{n-1}_{n-2})^{-1}(O)$. Along the $\BP^1$-fibers of $p^{n-1}_{n-2}$, this
line bundle is isomorphic to $\sO_{\BP^1}(3)$. Thus the difference from the equation $u^{(n)}=0$ of Proposition \ref{p.special}  is given by a cubic polynomial in $u^{(n-1)}$.
Thus it is of the form (\ref{e.ode}).

Conversely, the ODE (\ref{e.ode}) gives rise to a section of the line bundle  $T^{ p^{n-1}_{n-2}} \otimes \sO_{\BP \sJ^{n-2}}(1) $ on $(p^{n-1}_{n-2})^{-1}(O)$ for a suitable open subset  $O \subset J^{n-2}$, inducing a  splitting of the exact sequence. Thus it is associated with  a family of unbendable rational curves via Theorem \ref{t.Goursat}.\end{proof}

When $\dim M \leq 4$, unbendable rational curves of anti-canonical degree $3$ are essentially of Goursat type:

\begin{proposition}\label{p.dim4}
Let $M$ be a complex manifold of dimension $\leq 4$.
If $\sK$ is a bracket-generating family of unbendable rational curves of anti-canonical degree $3$, then it is of Goursat type. \end{proposition}

\begin{proof}
Under the assumption, the dimension of $\sK$ is equal to $\dim M \leq 4$ and the distribution $\sD$ is bracket-generating.
It is easy to see that a bracket-generating distribution of rank $2$ on a manifold of dimension $\leq 4$ is a Goursat distribution on its regular locus. Thus $\sK$ is of Goursat type. \end{proof}

\begin{example}\label{e.cubic}
Let $X \subset \BP^5$ be a nonsingular cubic hypersurface.
The family $\sK$ of lines on $X$ is a family of unbendable rational curves. The VMRT $\mathcal{C}_x$ at a general point $x\in X$ is irreducible and linearly nondegenerate in $\mathbb{P}T_xX$ (see Section 1.4.2 of \cite{Hw01}). The family $\mathcal{K}$ is bracket-generating by Proposition \ref{p.sT}. Thus it is of Goursat type by Proposition \ref{p.dim4} and there is an ODE of the type given in Theorem \ref{t.ode} associated with it.
However, it seems hard to find holomorphic functions $a_3, \ldots, a_0$ in Theorem \ref{t.ode} explicitly from the cubic equation defining $X$.
\end{example}

Now we discuss important features of the VMRT of unbendable rational curves of Goursat type.

\begin{proposition}\label{p.Goursat}
Use the notation of Definition \ref{d.vmrt} for a family of  unbendable rational curves of Goursat type on a manifold $M$.
Then the VMRT $\sC_{\alpha} \subset \BP T_{x} M$ for a general point $\alpha \in \sU$ is linearly nondegenerate  and the meromorphic distributions $\sT^k$ satisfy $[\sF, \sT^k] \subset \sT^k$ for all $k\geq 1$. \end{proposition}

\begin{proof}
We claim that ${\rm rank} (\sT^i) < {\rm rank}( \sT^{i+1})$ for $i \leq n - 2$ where $n$ is the dimension of $M$.  Since $\partial^k \sD$ has rank $k+2$ and by (\ref{e.dist})
$$\partial^k \sT^1 = \partial^k (\rho^{-1} \sD) = \rho^{-1}(\partial^k \sD)$$ has rank $k+3$ for each $0 \leq k \leq n-2$, the claim implies that  $\sT^{k+1} = \rho^{-1}(\partial^k \sD)$ for all $k \geq 0$, which proves the proposition.

To prove the claim, we may replace $\sK$ by any open subset in $\sK$. Thus by Theorem \ref{t.Goursat}, we may assume that the distribution $\sD$ on $\sK$ is biholomorphic to the distribution $\sJ^{n-2}$ on an open subset $O \subset J^{n-2}$.  Since  the derived distributions   of $\sJ^{n-2}$ are (see the last paragraph of Subsection \ref{ss.j})
$$\sJ^{n-2} \subset (p^{n-2}_{n-3})^{-1} \sJ^{n-3} \subset \cdots \subset (p^{n-2}_1)^{-1} \sJ^1 \subset TJ^{n-2},$$
we see that for a general element $v \in \sJ^{n-2}_y$,
$${\rm Levi}^{(p^{n-2}_i)^{-1} \sJ^i}(v, (p^{n-2}_i)^{-1} \sJ^i) \neq 0$$ for each $1 \leq i \leq n-2.$ By choosing a general $\alpha \in \rho^{-1}(y)$ such that ${\rm d}\rho(\sV_{\alpha}) $ corresponds to such a general element $v\in \sJ^{n-2}_y$, we see that ${\rm rank} (\sT^i) < {\rm rank}(\sT^{i+1})$ for $1 \leq i \leq n-2$ by induction.
\end{proof}

\begin{proposition}\label{p.GoursatIII}
Let $\sK$ be a family of unbendable rational curves of Goursat type on a complex manifold $M$ of dimension $n \geq 4$. Then for a general $y \in \sK$,
the third fundamental form of $\sC_{\alpha} \subset \BP T_{\mu(\alpha)} M$ vanishes for some $\alpha \in \rho^{-1}(y)$. \end{proposition}

\begin{proof}
Since $[\sF, \sT^2] \subset \sT^2$ and $[\sF, \sT^3] \subset \sT^3$ by Proposition \ref{p.Goursat}, the line bundle $\sT^3/\sT^2$ is trivial along $\rho^{-1}(y)$ for a general $y \in \sK$. Thus Proposition \ref{p.sT} (iii) implies that the relative normal spaces $\{N^{(3)}_{\sC_{\alpha}}, \alpha \in \rho^{-1}(y)\}$ define a vector bundle on $\rho^{-1}(y)$ isomorphic to  $$\sT_{3} \otimes T^* (\rho^{-1}(y)) \cong  \sO_{\BP^1}(-2).$$
Then as in Proposition \ref{p.FFvmrt} (iv),  the third fundamental forms of $\{\sC_{\alpha}, \alpha \in \rho^{-1}(y)\}$  is isomorphic to a section of the line bundle $$\Hom(\Sym^3 \sO_{\BP^1}(-1), \sO_{\BP^1}(-2)) = \sO_{\BP^1}(1).$$ It follows that the third fundamental form vanishes at some point of $\rho^{-1}(y)$. \end{proof}

\begin{lemma}\label{l.5}
Let $D$ be a bracket-generating distribution of rank 2 on a $5$-dimensional manifold $Y$. Then ${\rm rank}(\partial D) = 3$ and ${\rm rank} (\partial^{(2)} D) = {\rm rank}(\partial^2 D).$
\end{lemma}

\begin{proof}
We may assume that $D$ is regular at every point of $Y$.
Since ${\rm rank} (D) = 2,$ the image of ${\rm Levi}^D: \wedge^2 D \to TY/D$ has rank $1$ and $\partial D$ has rank 3.
If $\partial^{(2)} D$ has rank 5, then $\partial^2 D$ must have rank 5. It remains to exclude the case when ${\rm rank}(\partial^{(2)} D)= 4 $ and  ${\rm rank}(\partial^2 D) =5$. In this case, the homomorphism of vector bundles
$${\rm Levi}^{\partial D}: \wedge^2 (\partial D) \to TY/\partial D$$ is surjective on the regular locus of $D$.
For a general point $y \in Y$, choose a basis $u,v,w \in (\partial D)_y$ such that $u,v \in D_y$. Then ${\rm Levi}^{\partial D}(u,v) =0$, while ${\rm Levi}^{\partial D}(u,w)$ and $ {\rm Levi}^{\partial D}(v,w)$ have values in $\partial^{(2)}D/\partial D$. Thus ${\rm Levi}^{\partial D}$ cannot be surjective, a contradiction. \end{proof}

Now we have the following characterization of unbendable rational curves of Goursat type in dimension 5.

\begin{theorem}\label{t.test5}
Let $M$ be a 5-dimensional complex manifold and let $\sK$ be a bracket-generating family of  unbendable rational curves on $M$ with anti-canonical degree $3$. Then the following three statements are equivalent.
\begin{itemize}
\item[(1)] The growth vector of the distribution $\sD$ on  $\sK$ is $(2,3,4,5)$ (i.e., it is of Goursat type);
      \item[(2)] $[\sF, \sT^3] \subset \sT^3$;
      \item[(3)] for a general $y \in \sK$, the third fundamental form of $\sC_{\alpha}$ at $\tau( \alpha)$  vanishes for some $\alpha \in \rho^{-1}(y)$. \end{itemize} \end{theorem}

\begin{proof}
The implications $(1) \Rightarrow (2) \Rightarrow (3)$ are by Propositions \ref{p.Goursat} and \ref{p.GoursatIII}.

By Corollary \ref{c.T2},
$$\sT^2 = \partial^{(1)} \sT^1 = \rho^{-1}(\partial^{(1)} \sD).$$
Then (2) implies  $$\sT^3 = [\sT^1, \sT^2] + \sT^2 = \partial^{(2)} \sT^1 =  \rho^{-1}(\partial^{(2)} \sD),$$ where the first equality follows from Lemma \ref{l.ij} and the third equality follows from (\ref{e.dist}).  Thus $\partial^{(2)} \sD$ has rank $4$, which implies (1) by Lemma \ref{l.5}.

By Proposition \ref{p.FFvmrt} (iii), (3) implies that the image of the relative third fundamental forms of $\{\sC_{\alpha}, \alpha \in \rho^{-1}(y)\}$ span a line subbundle of $N_{C_y} \cong \sO_{\BP^1}(1) \oplus \sO_{\BP^1}^{\oplus 3}$ corresponding to one of the $\sO_{\BP^1}$-factors of $\sO_{\BP^1}^{\oplus 3}$. Thus we have a distinguished vector subbundle $N_{C_y}^{\sharp} \subset N_{C_y}$ isomorphic to $\sO_{\BP^1}(1) \oplus \sO_{\BP^1}^{\oplus 2}$ spanned by $\sO_{\BP^1}(1)$ and the images of the second and the third fundamental forms of $\sC_{\alpha}, \alpha \in \rho^{-1}(y)$. Define $$\sD^{\sharp}_y := H^0(C_y, N_{C_y}^{\sharp}) \subset \sD_y$$ for  general $y \in \sK$. This gives  a meromorphic distribution $\sD^{\sharp}$ on $\sK$ of rank $4$ containing $\partial \sD$ such that $\rho^{-1}\sD^{\sharp} =
\sT^3$ at general points of $\sU$. This implies (2).  \end{proof}

\section{Unbendable rational curves of Cartan type}\label{s.Cartan}

We introduce the following terminology for convenience.

\begin{definition}\label{d.Cartan}
Let $Y$ be a complex manifold of dimension 5.
A regular distribution $D \subset T Y$ of rank 2 with the growth vector $(2,3,5)$ is called a {\em Cartan distribution}. For a Cartan distribution $D \subset TY$,   the vector bundle homomorphisms given by Lie brackets
$$\wedge^2 D \to (\partial D)/D \mbox{ and }
D \otimes     (\partial D)/D \to TY/\partial D$$
are surjective. A family $\sK$ of unbendable rational curves in a 5-dimensional complex manifold is said to be {\em of Cartan type} if the distribution $\sD$ on $\sK$ in Definition \ref{d.sD} is a Cartan distribution on its regular locus.
\end{definition}

The following is straightforward.

\begin{lemma}\label{l.Cartan}
Let $D \subset TY$ be a bracket-generating distribution of rank 2 on a 5-dimensional complex manifold $Y$. Then the restriction of $D$ on its regular locus  is either a Cartan distribution or a Goursat distribution.
\end{lemma}

\begin{theorem}\label{t.test}
Let $\sK$ be a bracket-generating family of unbendable rational curves of anti-canonical degree 3 on a complex manifold of dimension 5.
Then it is of Cartan type if and only if for a general $y \in \sK$, the third fundamental form of $\sC_{\alpha}$ is nonzero at $\tau(\alpha)$ for each $\alpha \in \rho^{-1}(y)$. \end{theorem}

\begin{proof}
If the third fundamental form of $\sC_{\alpha}$ is nonzero at $\tau(\alpha)$ for each $\alpha \in \rho^{-1}(y)$, then $\sK$ is not of Goursat type by Theorem \ref{t.test5}. Thus it is of Cartan type by Lemma \ref{l.Cartan}.

Assuming that $\sK$ is of Cartan type, the distribution $\sD$ is a Cartan distribution in a neighborhood of a general point $y \in \sK$. Then $${\rm Levi}^{\partial \sD}(v, (\partial \sD)_y) \neq 0$$ for any nonzero  element $v \in \sD_y$ from Lemma \ref{l.5}. Since $\sT^2 = \rho^{-1} \partial \sD$ by Corollary \ref{c.T2}, this implies that the third fundamental form of $\sC_{\alpha}$ at $\tau(\alpha)$ is nonzero for each $\alpha \in \rho^{-1}(y)$ by Proposition \ref{p.sT} (ii) and Proposition \ref{p.rank2} (i). \end{proof}

There are many examples of Cartan distributions. In fact, a generic distribution of rank 2 on a 5-dimensional manifold is a Cartan distribution on its regular locus.
Applying Proposition \ref{p.rank2}, we obtain many examples of unbendable rational curves of Cartan type.
There is a particularly interesting class of examples related to contact structures.  To describe them, we need to recall some basics of contact structures.

\begin{definition}\label{d.contact}
For a complex manifold $M$, a regular distribution $H \subset TM$ of corank 1 is called a {\em contact structure} if ${\rm Ch}(H) =0$.  The line bundle $TM/H$ is called the {\em contact line bundle} and the line bundle-valued 1-form $\theta: TM \to TM/H$ is called the {\em contact form}.
Then ${\rm d} \theta |_H$ agrees with $-{\rm Levi}^H$.  When $n= \dim M$, there is a natural isomorphism $\wedge^n TM \cong (TM/H)^{(n+1)/2}.$ In particular, the dimension of $M$ should be an odd number. \end{definition}

\begin{example}\label{e.cotangent}
Let $Y$ be a complex manifold and let $\phi:  T^*Y \to Y$ be the projection of the cotangent bundle.
The natural 1-form $\theta^Y$ on $T^* Y$  sends a vector $v \in T_{\zeta}(T^*Y)$ at $\zeta \in T^*_y Y$ to  $$\theta^Y(v) = \zeta ( {\rm d}\phi (v)) \in \C.$$
Then ${\rm d} \theta^Y$ is the natural symplectic form on $T^*Y$. Let $\varphi: \BP T^*Y \to Y$ be the projective bundle and let $$\xi: T^*Y \setminus \mbox{(0-section)} \to \BP T^* Y$$ be the $\C^*$-bundle.
For a nonzero element $\zeta \in T_y^*Y$, denote by $[\zeta]\in \BP T^*_yY$  the image $\xi(\zeta)$   and by $\zeta^{\perp} \subset T_y Y$ the hyperplane annihilated by $\zeta.$ Define a hyperplane $H^Y_{[\zeta]} \subset T_{[\zeta]} (\BP T^*Y)$ by $$H^Y_{[\zeta]} = ({\rm d}_{[\zeta]} \varphi)^{-1} (\zeta^{\perp})$$ where ${\rm d}_{[\zeta]} \varphi: T_{[\zeta]}(\BP T^*Y) \to T_y Y$ is the differential of $\varphi$ at the point $[\zeta] \in \BP T^*_y Y$.
The resulting distribution $H^Y \subset T(\BP T^* Y)$ is a contact structure on $\BP
T^* Y$ satisfying ${\rm Ker}(\theta^Y) = \xi^{-1} H^Y$ on $T^*Y \setminus \mbox{(0-section)}$.  We denote the corresponding contact form by $$\vartheta^Y: T (\BP T^*Y) \to T(\BP T^*Y)/H^Y \cong \sO_{\BP T^*Y}(1).$$ \end{example}

We omit the proof of the following elementary lemma.

\begin{lemma}\label{l.reduction}
Let $S \subset X$ be a compact complex submanifold of a complex manifold $X$ and let $\vartheta: TX \to L$ be a surjective homomorphism to a line bundle $L$ on $X$ defining a distribution $\sH = {\rm Ker}(\vartheta) \subset TX$ of corank 1. Assume that \begin{itemize} \item[(1)] ${\rm Null}({\rm Levi}_x^{\sH}) := \{ v \in \sH_x, {\rm Levi}^{\sH}(v, u) = 0 \mbox{ for all } u \in \sH_x\}$ has the same dimension for all $x \in X$ and
  \item[(2)] ${\rm Null}({\rm Levi}_x^{\sH}) \cap T_x S =0$ for all $x \in S$. \end{itemize}
Then there exists a neighborhood $U \subset X$ of $S$ and a
submersion $\nu:U \to M$ to a complex manifold with a contact structure $H \subset TM$ such that $\sH|_{U}= \nu^{-1} H$ and $T^{\nu}_x = {\rm Null}({\rm Levi}_x^{\sH})$ for all $x \in U$. \end{lemma}

\begin{proposition}\label{p.contact}
Let $H \subset TM$ be a contact structure on a complex manifold of dimension 5.  Let $\sK$ be a bracket-generating family of unbendable rational curve of anti-canonical degree 3 whose members are tangent to $H$. Then $\sK$ is of Cartan type. \end{proposition}

\begin{proof}
The varieties of minimal rational tangents of $\sK$ at  general points of $M$ are contained in $\BP H \subset \BP TM$. Thus Proposition \ref{p.Goursat} implies that $\sK$ cannot be of Goursat type. Then it must be of Cartan type by Lemma \ref{l.Cartan}. \end{proof}

\begin{definition}\label{d.contactCartan}
A family of {\em contact unbendable rational curves of Cartan type} means  a family of unbendable rational curves of Cartan type described in Proposition \ref{p.contact}. \end{definition}

There are many examples of contact unbendable rational curves of Cartan type.
To see this, we need  the following properties of Cartan distributions  proved by Zelenko in \cite{Z99} and \cite{Z06}. They were  discussed also in Section 4.3 of \cite{BH} in a less precise form.

\begin{lemma}\label{l.Zelenko}
Let $D \subset TY$ be a Cartan distribution on a 5-dimensional complex manifold $Y$ and let $W \subset T^*Y$ be the subbundle of rank $2$ annihilating $\partial D$. Using the notation of  Example \ref{e.cotangent}, let $\sigma:= {\rm d}\theta^Y|_W$ be the 2-form on $W$ given by  the restriction of ${\rm d} \theta^{Y}$. Then the null-space $${\rm Null}(\sigma)_w := \{ u \in T_w W, \sigma(u, v) = 0 \mbox{ for all } v \in T_w W\} $$ is 1-dimensional for each nonzero  $w \in W,$ defining a line subbundle $\sV^{\sharp} \subset TW$ outside the $0$-section which satisfies the following.
\begin{itemize}
\item[(i)] $\theta^Y(\sV^{\sharp}) =0.$ \item[(ii)] $\sV^{\sharp} \cap T^{\eta} = 0$ where $\eta: W \to Y$ is the natural projection.
    \item[(iii)] $\sV^{\sharp} \subset \eta^{-1} D.$
      \end{itemize} \end{lemma}

\begin{proof}
All the statements are proved in \cite{Z06} (Zelenko worked in the setting of real differentiable manifolds, but the arguments are valid also in complex analytic setting).
That  ${\rm Null}(\sigma)_w$ is a 1-dimensional subspace in ${\rm Ker}(\theta^Y)$  is the equation (3.13) of \cite{Z06} which was proved in Proposition 2.2 and Corollary 2.1 of \cite{Z99}. (ii) and (iii) follow from the equations  (3.13) and   (3.16) in the proof of Proposition 3.1 in  \cite{Z06}.
\end{proof}

The following theorem says that each germ of a Cartan distribution determines in a natural way a germ of a family of contact unbendable rational curves  of Cartan type.

\begin{theorem}\label{t.Zelenko}
Let $D \subset T Y$ be a Cartan distribution. Let $\vartheta^Y$ be the canonical contact form on $\BP T^* Y$ from Example \ref{e.cotangent}. Let $\varrho: \BP W \to Y$ be the projectivization of the vector bundle $W \subset T^*Y$ in Lemma \ref{l.Zelenko} and let $\vartheta : T\BP W \to \sO_{\BP W}(1)$ be the line-bundle valued 1-form on $\BP W \subset \BP T^*Y$ obtained from the restriction of $\vartheta^Y$. Then the following holds.
\begin{itemize}
\item[(1)] ${\rm Ker}(\vartheta)\subset T \BP W$ is a distribution of corank 1 on $\BP W$.
    \item[(2)] In the terminology of Lemma \ref{l.reduction}, define $$ \sV^{\flat}_{x} := {\rm Null}({\rm Levi}_{x}^{{\rm Ker}(\vartheta)}) \mbox{ for } x \in \BP W.$$ Then $\dim \sV^{\flat}_{x} =1 $ for all $x \in\BP W$ defining a line subbundle  $\sV^{\flat} \subset T \BP W$, which is  transversal to $\sF^{\flat} := T^{\varrho}$ and satisfies $\sV^{\flat} \subset \varrho^{-1} D$.
        \item[(3)]
For any point $y_o \in Y$, there exists a neighborhood $\sK^{\flat} \subset Y$ of $y_o$ such that
 the leaves of $\sV^{\flat}$ define a submersion $\nu: \sU^{\flat}:= \varrho^{-1}(\sK^{\flat}) \to M^{\flat}$ to a 5-dimensional complex manifold $M^{\flat}$ equipped with a contact structure $H^{\flat}$ satisfying $\nu^{-1} H^{\flat} = {\rm Ker}(\vartheta)|_{\sU^{\flat}}$.
 \item[(4)] In (3), for any $y \in \sK^{\flat}$, the morphism $\nu$ sends $\varrho^{-1}(y)$ to an unbendable rational curve $C^{\flat}_y \subset M^{\flat}$ of anti-canonical degree $3$ satisfying $T C^{\flat}_y \subset H^{\flat}|_{C^{\flat}_y}$.
 \item[(5)] Regarding $$ \sK^{\flat} \stackrel{\varrho}{\longleftarrow} \sU^{\flat} \stackrel{\nu}{\longrightarrow} M^{\flat}$$ in (3) as a family of unbendable rational curves on $M^{\flat}$ via (4), we have $D|_{\sK^{\flat}} = \sD$ where $\sD$ is the distribution determined by the family of unbendable rational curves on $M^{\flat}$ as in Definition \ref{d.sD}.
\end{itemize}  \end{theorem}

\begin{proof}
(1) is immediate from the definition of $\theta^Y$ and the submersion $\varrho: \BP W \to Y$.

(2) is a reformulation of Lemma  \ref{l.Zelenko}: the line bundle $\sV^{\flat}$ is just the image of $\sV^{\sharp}$ of Lemma \ref{l.Zelenko} under the projection $\lambda: W \setminus \mbox{(0-section)} \to \BP W$. To see this, let $\theta$ be the restriction of $\theta^Y$ to $W$ such that
${\rm Ker}(\theta) = \lambda^{-1}({\rm Ker}(\vartheta))$ from ${\rm Ker}(\theta^Y) = \xi^{-1}{\rm Ker}(\vartheta^Y)$ in Example \ref{e.cotangent}.  For a nonzero $w \in W$, it is easy to see that $${\rm Null}({\rm Levi}_w^{{\rm Ker}(\theta)}) = {\rm Null}({\rm d}\theta|_{{\rm Ker}(\theta)_w}).$$
Thus it suffices to show  \begin{eqnarray}\label{e.null} {\rm Null}({\rm d} \theta|_{{\rm Ker}(\theta)_w}) & \subset & {\rm Null}({\rm d} \theta)_w + T_w^{\lambda}.\end{eqnarray} Fix $e \in T_w^{\lambda}$ and $u \in T_wW$ satisfying ${\rm d} \theta (e, u) =1$ and $\theta(u) \neq 0.$
For any $v \in {\rm Null}({\rm d} \theta|_{{\rm Ker}(\theta)_w})$, we have
\begin{equation}\label{e.theta} {\rm d} \theta (v - {\rm d} \theta (v, u) e, u ) = 0. \end{equation}  Note that $v - {\rm d} \theta (v, u) e \in {\rm Null}({\rm Levi}_w^{{\rm Ker}(\theta)})$ from $T^{\lambda}_w \subset {\rm Null}({\rm Levi}_w^{{\rm Ker}(\theta)}).$  Thus (\ref{e.theta})
 implies that $v - {\rm d} \theta(v, u) e \in {\rm Null}({\rm d} \theta)_w$. This proves (\ref{e.null}) and (2).

(3) is a consequence of (2) and Lemma \ref{l.reduction}.

In (4), the inclusion $T C^{\flat}_y \subset H|_{C^{\flat}_y}$ is immediate from $T^{\varrho} \subset {\rm Ker}(\vartheta)$.
 The  normal bundle $N_{C^{\flat}_y}$ of $C^{\flat}_y$ in $M^{\flat}$ is a quotient bundle of the normal bundle $N_{\varrho^{-1}(y)}$ of $\varrho^{-1}(y)$ in $\sU^{\flat}$ which is a trivial vector bundle. Thus $N_{C^{\flat}_y}$ is semi-positive. From $$\wedge^5 T M^{\flat} \cong (TM/H)^{\otimes 3} $$ in Definition \ref{d.contact} and $$(TM/H)|_{C^{\flat}_y} \cong \sO_{\BP^1}(1)$$ in Example \ref{e.cotangent},  we have $\wedge^5 T M^{\flat}|_{C^{\flat}_y} \cong \sO_{\BP^1}(3)$, which implies
  $$N_{C^{\flat}_y} \cong \sO_{\BP^1}(1) \oplus \sO_{\BP^1}^{\oplus 3}.$$
This shows that $C^{\flat}_y$ is an unbendable rational curve.

To see (5), note that  $\sV^{\flat} + \sF^{\flat} \subset \varrho^{-1}D$ from (2). Combining it with Proposition \ref{p.HM}, we have
$\sV^{\flat} + \sF^{\flat} \subset \varrho^{-1}(D \cap \sD)$.   This inclusion must be an identity if ${\rm rank}(D \cap \sD) \leq 1$.
But then $\sF^{\flat}$ must be the Cauchy characteristic of $\sV^{\flat} +\sF^{\flat}$, a contradiction to Proposition \ref{p.sT} (i). It follows that ${\rm rank}(D \cap \sD) =2,$ implying $D = \sD$.
\end{proof}

\begin{remark}
The construction in Theorem \ref{t.Zelenko} is inspired by the
idea of Jacobi curves in \cite{Z99} and \cite{Z06}. Our argument is a translation of Zelenko's  symplectic approach into   contact geometry. Zelenko's Jacobi curves can be interpreted as the images of VMRT under the Gauss map. \end{remark}

The next theorem says that the  unbendable rational curves of Cartan type arising from the construction of Theorem \ref{t.Zelenko} cover all examples of contact unbendable rational curves of Cartan type.

\begin{theorem}\label{t.model}
In the setting of Proposition \ref{p.contact}, let $Y \subset \sK$ be the regular locus of $\sD$ and let $D \subset TY$ be the Cartan distribution given by the restriction of $\sD$ to $Y$. We can apply Theorem \ref{t.Zelenko} to $D$ to obtain  $\sK^{\flat}, \sU^{\flat}, M^{\flat}$ and $H^{\flat}$. Then there exists an open embedding $\chi: \sU^{\flat} \subset \sU$ satisfying $$ \varrho = \rho \circ \chi \mbox{ and } {\rm d} \chi( {\rm Ker}(\vartheta)) = \mu^{-1}H |_{\chi(\sU^{\flat})}.$$ In particular, there is an open embedding $M^{\flat} \subset M$ satisfying $H^{\flat} = H|_{M^{\flat}}$ and a commutative diagram $$ \begin{array}{ccccc} \sK^{\flat} &\stackrel{\varrho}{\leftarrow} & \sU^{\flat} & \stackrel{\nu}{\rightarrow} & M^{\flat} \\
\cap & & \downarrow \chi & & \cap \\ \sK & \stackrel{\rho}{\leftarrow} & \sU & \stackrel{\mu}{\rightarrow} & M.\end{array}$$
\end{theorem}

\begin{proof}
To prove the theorem, we may replace $\sK$ by $Y$ and assume that $D = \sD$ is a Cartan distribution.
Then $\sT^2$ is a distribution on $\sU$ satisfying $\sT^2 = \rho^{-1}\partial D$ by Corollary \ref{c.T2}.

 Since elements of $\sK$ are tangent to $H$, the germ $\sC_{\alpha}$ is contained in $\BP H_x $ for all $x \in M$ and $\alpha \in \mu^{-1}(x)$. Then $\sT^m \subset \mu^{-1}H$ for all $m \geq 0$ by Proposition \ref{p.sT} (ii). In particular, the annihilators of $\mu^{-1}H$ define a section $\sigma: \sU \to \BP (T\sU/\sT^2)^*$ of the projection $\BP (T\sU/
\sT^2)^* \to \sU$. Let $\Sigma \subset \BP(T\sU/\sT^2)^*$ be the image of $\sigma$.  It is easy to check from the definition of $\vartheta^{\sU}$ that \begin{eqnarray}\label{e.H} {\rm Ker}(\vartheta^{\sU}|_{\Sigma}) &\cong& \mu^{-1}H \end{eqnarray} under the biholomorphic map $\Sigma \cong \sU$.

Let $W \subset T^*Y$ be the annihilator of $\partial D$ and define a holomorphic map $\eta: \sU \to \BP W $ by sending $\alpha \in \sU$ to the annihilator of the hyperplane ${\rm d} \rho (\mu^{-1}H)_{\alpha} \subset T_y Y$ containing $(\partial D)_y $ with $y = \rho(\alpha)$. From the commutative diagram of $\BP^1$-bundles
$$ \begin{array}{ccc}
\sU & \stackrel{\eta}{\to} & \BP W \\
\downarrow \rho & & \downarrow \\
Y & = & Y, \end{array}$$
if $\eta$ is not biholomorphic, it  contracts  each fiber of $\rho$ to  a single point. In the latter case, we obtain a distribution $H'$ on $Y$ such that $\mu^{-1}H = \rho^{-1}H'$, which implies $\sF \subset {\rm Ch}(\mu^{-1}H)$, a contradiction to
$${\rm Ch}(\mu^{-1}H) = \mu^{-1}{\rm Ch}(H)= T^{\mu}.$$
  It follows that $\eta$ is a biholomorphic map.

The homomorphism  of vector bundles $\rho^* (T^*Y) \to T^* \sU$ induces an isomorphism of vector bundles $\rho^* W \cong (T\sU/\sT^2)^*$ because $\sT^2 = \rho^{-1}\partial D$. Let $\xi: \BP(T \sU/\sT^2)^* \to \BP W$ be the composition of  the holomorphic maps $$\BP(T \sU/\sT^2)^* \stackrel{\cong}{\rightarrow} \rho^* \BP W \to \BP W.$$ It is easy to see that \begin{eqnarray}\label{e.xi} \xi^{-1}{\rm Ker}(\vartheta^Y|_{\BP W}) &=& {\rm Ker}(\vartheta^{\sU}|_{\BP (T \sU/\sT^2)^*})\end{eqnarray} and
$\xi \circ \sigma = \eta.$ Equations (\ref{e.H}) and (\ref{e.xi}) say that the biholomorphism $\eta$ sends $\mu^{-1}H$ to ${\rm Ker}(\vartheta^Y|_{\BP W})$. By the definitions of $H^{\flat}$ and $M^{\flat}$ in Theorem \ref{t.Zelenko}, this implies that the inclusions $\sU^{\flat} \subset \sU$ and $M^{\flat} \subset M$ obtained from the inverse of $\eta$ satisfy the required properties. \end{proof}

\section{Lines on quartic 5-folds}\label{s.quartic}

\subsection{Main result and the strategy of its proof}

The main result of this section is the following.

\begin{theorem}\label{t.general}
A general hypersurface $X \subset\BP^6$ of degree $4$ and a general projective line $\ell$ on $X$ have the following properties.
\begin{itemize}
\item[(1)] $X$ is smooth along $\ell$.
\item[(2)]  The normal bundle $N_{\ell}$ of $\ell \subset X$ is isomorphic to  $\sO(1)\oplus\sO^{\oplus 3}$.
\item[(3)] The second fundamental form ${\rm FF}^2_{\ell_x, \sC_x}$ is nonzero for all $x\in \ell$. \item[(4)] The third fundamental form ${\rm FF}^3_{\ell_x, \sC_x}$ is nonzero for all $x\in \ell$.
    \item[(5)] The fourth fundamental form  ${\rm FF}^4_{\ell_P, \sC_P}$ is nonzero for some $P \in \ell$ and consequently, the VMRT $\sC_P$ at $P$ is linearly nondegenerate in $\BP T_P X.$
\end{itemize}  In particular, general lines on $X$ are unbendable rational curves of Cartan type by Proposition \ref{p.sT} and Theorem \ref{t.test}.
\end{theorem}

To prove  Theorem \ref{t.general}, it suffices to prove the next theorem because the properties (1)--(5) are  open conditions in the space of pairs $(X, \ell)$ consisting of a  quartic hypersurface  $X \subset \BP^6$ and an unbendable line $\ell$ contained in the smooth locus of $ X.$

\begin{theorem}\label{t.special}
Let $\BP^6$ be the projective space with homogeneous coordinates $x_0,\ldots,x_6$. Define
\begin{eqnarray*}
f& :=&\sum_{i=2}^6x_i^4
+x_0^3x_2+x_1^3x_3+x_0^2x_1x_4+x_0x_1^2x_5 +x_0^2x_2^2+x_1^2x_3^2 \\
&&+(x_0^2+x_0x_1+x_1^2)x_6^2
+x_0x_4^3+x_1x_5^3+(x_0+x_1)x_6^3. \nonumber
\end{eqnarray*}
Then $\ell:=\cap_{i=2}^6\zero(x_i)$ is a projective line on the hypersurface $X:=\zero(f)$ of degree $4$ with the following properties.
\begin{itemize}
\item[(1)] $X$ is smooth along $\ell$.
\item[(2)]  The normal bundle $N_{\ell}$ of $\ell \subset X$ is isomorphic to  $\sO(1)\oplus\sO^{\oplus 3}$.
\item[(3)] The second fundamental form ${\rm FF}^2_{\ell_x, \sC_x}$ is nonzero for all $x\in \ell$. \item[(4)] The third fundamental form ${\rm FF}^3_{\ell_x, \sC_x}$ is nonzero for all $x\in \ell$.
    \item[(5)] The fourth fundamental form  ${\rm FF}^4_{\ell_P, \sC_P}$ is nonzero for some $P \in \ell$ and consequently, the VMRT $\sC_P$ at $P$ is linearly nondegenerate in $\BP T_P X.$
\end{itemize}
\end{theorem}

\begin{lemma}\label{l.smooth}
The hypersurface $X$ is smooth along the projective line $\ell$.
\end{lemma}

\begin{proof}
Take any point $x=[t_0:t_1:0:\cdots:0]\in \ell$. Then the tangent space of the affine cone $\hat{X} \subset \C^7$ of $X$ at $x$ is
\begin{eqnarray}\label{eqn-tantent}
\{t_0^3x_2+t_1^3x_3+t_0^2t_1x_4+t_0t_1^2x_5=0\},
\end{eqnarray}
which is properly contained in $\C^7$. Thus $X$ is smooth at $x$.
\end{proof}

\begin{notation}\label{n.vmrt} We can identify the space of lines through a point $x \in \BP^n$  with $\BP T_x \BP^n$ in a natural way.  Using this identification and abusing the notation, let us denote by $\sC_x \subset \BP T_x X \subset \BP T_x \BP^n$ the subscheme of lines lying on $X$ passing through $x$. If  $\ell$ is an unbendable rational curve, the germ of $\sC_x$ at $\ell_x$ agrees with the germ $\sC_\alpha$ of the VMRT at $\alpha=\mathbb{P}T_x\ell$ as defined in Definition \ref{d.vmrt}. \end{notation}

The next lemma is easy to check.

\begin{lemma}\label{l.involution}
Define an involution $\gamma$ of $\BP^6$ by $$[\lambda_0: \lambda_1: \lambda_2: \lambda_3: \lambda_4: \lambda_5:\lambda_6] \mapsto [\lambda_1: \lambda_0: \lambda_3: \lambda_2: \lambda_5: \lambda_4: \lambda_6].$$ Then we have $\gamma\cdot [\C f]= [\C f]$, $\gamma\cdot X=X$, $\gamma\cdot \ell =\ell$, and the points of $\ell$ fixed by $\gamma$  are $$P:=[1:1:0:\cdots:0] \mbox{ and } Q:=[-1:1:0:\cdots:0].$$
\end{lemma}

In the terminology of Notation \ref{n.vmrt} and Lemma \ref{l.involution}, we state the following two propositions, whose proofs are to be given later.

\begin{proposition}\label{p.p_0}
The subscheme $\sC_P \subset \BP T_P X$ is smooth at $\ell_P$, the germ of $\sC_P$ at $\ell_P$ is of dimension 1 and its fundamental forms at $\ell_P$,  $${\rm FF}^2_{\ell_P, \sC_P}, \  {\rm FF}^3_{\ell_P, \sC_P}, \mbox{ and } {\rm FF}^4_{\ell_P, \sC_P}$$ are nonzero.
\end{proposition}

\begin{proposition}\label{p.p_1} The subscheme $\sC_Q \subset \BP T_Q X$ is smooth at $\ell_Q$, the germ of $\sC_Q$ at $\ell_Q$ is of dimension 1  and
its fundamental forms at $\ell_Q$,  $${\rm FF}^2_{\ell_Q, \sC_Q}  \mbox{ and } {\rm FF}^3_{\ell_Q, \sC_Q}$$ are nonzero.
\end{proposition}

 Theorem \ref{t.special} can be proved from these two propositions as follows.

 \begin{proof}[Proof of Theorem \ref{t.special}]
 (1) is Lemma \ref{l.smooth}.

 Since the normal bundle of $\ell$ in $\BP^6$ is isomorphic to $\sO_{\ell}(1)^{\oplus 5}$, its normal bundle  $N_{\ell}$ in $X$ is isomorphic to $\oplus_{i=1}^4\sO(d_i)$ with each $d_i\leq 1$. By the adjunction formula, the anti-canonical degree of $\ell$ in $X$ is $\deg T\BP^6|_{\ell} -\deg X=3$ which implies $\deg N_{\ell}=\sum_{i=1}^4d_i=1$.
By the deformation theory of rational curves, there is a natural identification
$$T_{\ell_P}\sC_P=H^0(\ell, N_{\ell}\otimes {\bf m}_P))=H^0(\BP^1, \oplus_{i=1}^4\sO(d_i-1)).$$
By Proposition \ref{p.p_0}, we have $\dim T_{\ell_P}\sC_P=1$. By reordering $d_i$, we deduce from the above formula  that $d_1=1$ and $d_i\leq 0$ for $i\neq 1$. Since $\sum_{i=1}^4d_i=1$, we have $N_{\ell} \cong \sO(1)\oplus\sO^{\oplus 3}$. This proves (2).

 (3) follows from Proposition \ref{p.p_0} because the second fundamental form of the VMRT is unchanged along an unbendable rational curve by Proposition \ref{p.FFvmrt} (ii).

 To check (4), we use Proposition \ref{p.FFvmrt} (iii) to obtain a holomorphic line subbundle $L$ in the $\sO(1)^{\oplus 3}$-factor of $N_{\ell}$ in (3) such that the third fundamental form ${\rm FF}^3_{\ell_x, \sC_x}, x \in \ell$ is an element of $ H^0(\ell, L)$. In particular, either ${\rm FF}^3_{\ell_x, \sC_x}$ is zero for all $x \in \ell$ or there is at most one point $z\in \ell$ such that ${\rm FF}^3_{\ell_z, \sC_z} =0$. By Propositions \ref{p.p_0} and \ref{p.p_1}, we have ${\rm FF}^3_{\ell_z, \sC_z} =0$ for at most one $z\in \ell\setminus\{P, Q\}$. But if such a point $z$ exists, as the involution $\gamma$  in Lemma \ref{l.involution} induces an isomorphism between ${\rm FF}^3_{\ell_x, \sC_x}$ with ${\rm FF}^3_{\ell_{\gamma(x)}, \sC_{\gamma(x)}}$ for all $x\in \ell$, we have $\gamma(z)\neq z$ and $${\rm FF}^3_{\ell_{\gamma(z)}, \sC_{\gamma(z)}} =0,$$ a contradiction. This proves (4).

  Finally, (5) is from Proposition \ref{p.p_0}. \end{proof}

\subsection{Proof of Proposition \ref{p.p_0}}

\begin{lemma}\label{l.p_0Z}
Let $\BP^5 \subset \BP^6$ be the hyperplane $\zero(x_0)$ with homogenous coordinates $[y_1: \cdots : y_6]$ given by $y_i=x_i, 1\leq i \leq 6.$
Let $\psi: \BP^6\dashrightarrow \BP^5$ be the projection with vertex $P$ which induces an isomorphism $[\psi]:\BP(T_{P}\BP^6)\rightarrow\BP^5$.
Then $[\psi]$ identifies $$\sC_{P} \subset \BP T_{P}X \subset \BP (T_{P} \BP^6)$$ with the subscheme  $Z\subset\BP^5$ defined by the system of homogeneous equations $h_1=h_2=h_3=h_4=0$, where
\begin{eqnarray*}
 && h_1:= y_2+y_3+y_4+y_5, \\
 && h_2:=3y_1y_3+y_1y_4+2y_1y_5+y_2^2+y_3^2+3y_6^2,  \\
 && h_3:=3y_1^2y_3+y_1^2y_5+2y_1y_3^2+3y_1y_6^2+y_4^3+y_5^3+2y_6^3,\\
 && h_4:=\sum_{i=2}^6y_i^4+y_1^3y_3+y_1^2y_3^2+y_1^2y_6^2+y_1y_5^3+y_1y_6^3.
\end{eqnarray*}
Under this identification, the point $\ell_P \in \sC_{P}$ corresponds to $\alpha:=[1:0:\cdots:0] \in Z.$
\end{lemma}

\begin{proof}
The projective line  joining $P$
and $[y_1:\cdots:y_5] \in \BP^5$ intersects $\BP^6\setminus \BP^5$ along
$$\{[1: \lambda y_1+1 : \lambda y_2: \cdots : \lambda y_6], \lambda \in \C\}.$$ This line lies on $X$, i.e., belongs to $\sC_{P}$ if and only if $$f(1, \lambda y_1+1, \lambda y_2,\ldots, \lambda y_6)=\sum_{k=1}^4h_k\lambda^k=0$$ for all $\lambda\in\C$.  It follows that  $h_1=h_2=h_3=h_4=0$ describes the image of the subscheme $\sC_{P}$ under the identification $[\psi]$. It is clear that $\ell_P$ is sent to $\alpha$.
\end{proof}

The following translation of Lemma \ref{l.p_0Z} into affine coordinates is immediate.

\begin{lemma}\label{l.p_0affine}
Let $\A^5 \subset \BP^5$ be the complement of  $\zero(y_1)$ with the affine coordinates $z_i:=\frac{y_i}{y_1}$, $i=2,\ldots, 6$. Then $\alpha \in Z$ is defined by $z_2= \cdots = z_6 =0$ and the subscheme $Z\cap\A^5$  is defined by $g_1=g_2=g_3=g_4=0$, where
\begin{eqnarray*}
 && g_1:= z_2+z_3+z_4+z_5, \nonumber \\
 && g_2:=3z_3+z_4+2z_5+z_2^2+z_3^2+3z_6^2,  \label{eqn-Z-intersect-A5}\\
 && g_3:=3z_3+z_5+2z_3^2+3z_6^2+z_4^3+z_5^3+2z_6^3, \nonumber\\
 && g_4:=\sum_{i=2}^6z_i^4+z_3+z_3^2+z_6^2+z_5^3+z_6^3. \nonumber
\end{eqnarray*} \end{lemma}

\begin{lemma}\label{l.p_0tangent}
In Lemma \ref{l.p_0affine}, the Zariski tangent space  $T_{\alpha}Z \cong T_{\ell_P} \sC_{P}$ corresponds to the 1-dimensional subspace $z_2= z_3=z_4 =z_5=0$.  In particular, the subscheme $\sC_P$ is of dimension 1 and smooth at $\ell_P$.
\end{lemma}

\begin{proof}
This is immediate from
\begin{eqnarray*}
&& T_\alpha(\zero(g_1))=\{z_2+z_3+z_4+z_5=0\},\\
&& T_\alpha(\zero(g_2))=\{3z_3+z_4+2z_5=0\}, \\
&& T_\alpha(\zero(g_3))=\{3z_3+z_5=0\}, \\
&& T_\alpha(\zero(z_4))=\{z_3=0\}.
\end{eqnarray*}
\end{proof}

\begin{lemma}\label{l.p_0abcd}
By Lemma \ref{l.p_0tangent}, the function $t:= z_6|_{Z \cap \A^5}$ gives a local coordinate on a neighborhood $U \subset Z$ of $\alpha$.
Write \begin{eqnarray}\label{e.taylor}
z_i &=& a_it+b_it^2+c_it^3+d_it^4+O(t^5)
\end{eqnarray} on  $U$ for all $ i=2,\ldots,6.$ Then \begin{eqnarray*}
\mathbf{a}=(a_2,a_3,a_4,a_5, a_6), & \mathbf{b}=(b_2, b_3,b_4,b_5,b_6), \\
\mathbf{c}=(c_2, c_3, c_4,c_5, c_6), & \mathbf{d}=(d_2,d_3,d_4,d_5,d_6). \nonumber
\end{eqnarray*} are given by
\begin{eqnarray}\label{e.abcd}
\mathbf{a}=(0, 0, 0, 0, 1), & \mathbf{b}=(1,-1, 0, 0, 0), \\
\mathbf{c}=(-1, -1, 1, 1, 0), & \mathbf{d}=(2, -2, -4, 4, 0). \nonumber
\end{eqnarray}
 In particular, the vectors $\mathbf{a}, \mathbf{b}$, $\mathbf{c}$ and $\mathbf{d}$ are linearly independent.
\end{lemma}

\begin{proof}
From $t =z_6|_U$, we already have $a_6=1$ and $b_6=c_6=d_6=0$.
Recall that $g_k(z_2,\ldots, z_6)=0$ on $Z\cap\A^5$ for $k=1,\ldots,4$. Putting (\ref{e.taylor}) in these equations, we obtain $g_k(z_2,\ldots, z_6)=\sum_{i=1}^4 G_{k, i}t^i+O(t^5)$ on $U$, where each $G_{k, i}$ is an explicit polynomials of $a_2,b_2, c_2, d_2,\ldots, a_6,b_6,c_6, d_6$.
 Since these coefficients $G_{k, i}$ of $t^i$ vanish, we can determine the vectors $\mathbf{a}$, $\mathbf{b}$, $\mathbf{c}$ and $\mathbf{d}$ one by one to obtain (\ref{e.abcd}).
\end{proof}

\begin{lemma}\label{l.0f}
The second, the third and the fourth fundamental forms of $Z$ at $\alpha$ are nonzero.
\end{lemma}

\begin{proof}
By the definition of vectors $\mathbf{a}$, $\mathbf{b}$, $\mathbf{c}$ and $\mathbf{d}$, we have
\begin{eqnarray*}
T_\alpha(Z\cap\A^5) &=& \C\mathbf{a}, \\
 T^{(2)}_\alpha(Z\cap\A^5)&=&\C\mathbf{a}+\C\mathbf{b}, \\
 T^{(3)}_\alpha(Z\cap\A^5)&=& \C\mathbf{a}+\C\mathbf{b}+\C\mathbf{c}, \\
 T^{(4)}_\alpha(Z\cap\A^5)&=&\C\mathbf{a}+\C\mathbf{b}+\C\mathbf{c}+\C\mathbf{d},
\end{eqnarray*}
where $T^{(k)}_\alpha(Z\cap\A^5)$ is the $k$-th osculating space of $Z\cap\A^5$ at $\alpha$. By Lemma \ref{l.p_0abcd}, we obtain the result.
\end{proof}

Lemma \ref{l.p_0tangent} and Lemma \ref{l.0f} complete the proof of Proposition \ref{p.p_0}.

\subsection{Proof of Proposition \ref{p.p_1}}

The proof of Proposition \ref{p.p_1} is completely parallel to the proof of Proposition \ref{p.p_0}. It consists of the following lemmata. We skip the proofs which are direct translations of those in   the proof of Proposition \ref{p.p_0}.

\begin{lemma}\label{l.p_1Z}
Let $\BP^5 \subset \BP^6$ be the hyperplane $\zero(x_0)$ with homogenous coordinates $[y_1: \cdots : y_6]$ given by $y_i=x_i, 1\leq i \leq 6.$
Let $\psi: \BP^6\dashrightarrow \BP^5$ be the projection with vertex $Q$ which induces an isomorphism $[\psi]:\BP(T_{Q}\BP^6)\rightarrow\BP^5$.
Then $[\psi]$ identifies $$\sC_{Q} \subset \BP T_{Q}X \subset \BP (T_{Q} \BP^6)$$ with the subscheme  $Z\subset\BP^5$ defined by the system of homogeneous equations $h_1=h_2=h_3=h_4=0$, where
\begin{eqnarray*}
 && h_1:= -y_2+y_3+y_4-y_5, \\
 && h_2:=3y_1y_3+y_1y_4-2y_1y_5+y_2^2+y_3^2+y_6^2,  \\
 && h_3:=3y_1^2y_3-y_1^2y_5+2y_1y_3^2+y_1y_6^2-y_4^3+y_5^3,\\
 && h_4:=\sum_{i=2}^6y_i^4+y_1^3y_3+y_1^2y_3^2+y_1^2y_6^2+y_1y_5^3+y_1y_6^3.
\end{eqnarray*}
Under this identification, the line $\ell_Q \in \sC_{Q}$ corresponds to $\alpha:=[1:0:\cdots:0] \in Z.$
\end{lemma}

\begin{lemma}\label{l.p_1affine}
Let $\A^5 \subset \BP^5$ be the complement of  $\zero(y_1)$ with the affine coordinates $z_i:=\frac{y_i}{y_1}$, $i=2,\ldots, 6$. Then $\alpha \in Z$ is defined by $z_2= \cdots = z_6 =0$ and the subscheme $Z\cap\A^5$  is defined by $g_1=g_2=g_3=g_4=0$, where
\begin{eqnarray*}
 && g_1:= -z_2+z_3+z_4-z_5, \nonumber \\
 && g_2:=3z_3+z_4-2z_5+z_2^2+z_3^2+z_6^2,  \label{eqn-Z-intersect-A5-p2}\\
 && g_3:=3z_3-z_5+2z_3^2+z_6^2-z_4^3+z_5^3, \nonumber\\
 && g_4:=\sum_{i=2}^6z_i^4+z_3+z_3^2+z_6^2+z_5^3+z_6^3. \nonumber
\end{eqnarray*} \end{lemma}

\begin{lemma}\label{l.p_1tangent}
In Lemma \ref{l.p_1affine}, the Zariski tangent space  $T_{\alpha}Z \cong T_{\ell_Q} \sC_{Q}$ corresponds to the 1-dimensional subspace $z_2= z_3=z_4 =z_5=0$.  In particular, the subscheme $\sC_{Q}$ is of dimension 1 and  smooth at $\ell_Q$.
\end{lemma}

\begin{lemma}\label{l.p_1abcd}
By Lemma \ref{l.p_1tangent}, the function $t:= z_6|_{Z \cap \A^5}$ gives a local coordinate on a neighborhood $U \subset Z$ of $\alpha$.
Write \begin{eqnarray*}
z_i&=&a_it+b_it^2+c_it^3+d_it^4+O(t^5)
\end{eqnarray*} on  $U$ for all $ i=2,\ldots,6.$ Then \begin{eqnarray*}
\mathbf{a}=(a_2,a_3,a_4,a_5, a_6), & \mathbf{b}=(b_2, b_3,b_4,b_5,b_6), \\
\mathbf{c}=(c_2, c_3, c_4,c_5, c_6), & \mathbf{d}=(d_2,d_3,d_4,d_5,d_6). \nonumber
\end{eqnarray*} are given by
\begin{eqnarray*}
\mathbf{a}=(0, 0, 0, 0, 1), & \mathbf{b}=(-1,-1, -2, -2, 0), \\
\mathbf{c}=(-1, -1, -3, -3, 0), & \mathbf{d}=(-2, -2, -4, -4, 0). \nonumber
\end{eqnarray*}
In particular, the vectors $\mathbf{a}, \mathbf{b}$ and $\mathbf{c}$ are linearly independent, and we have $\mathbf{d}=2\mathbf{b}$.
\end{lemma}

\begin{lemma}\label{l.1f}
The second and the third  fundamental forms of $Z$ at $\alpha$ are nonzero.
\end{lemma}

Lemma \ref{l.p_1tangent} and Lemma \ref{l.1f} complete the proof of Proposition \ref{p.p_1}.

\bigskip
Center for Complex Geometry,
Institute for Basic Science (IBS),
Daejeon 34126, Republic of Korea

\medskip
jmhwang@ibs.re.kr

qifengli@ibs.re.kr
 \end{document}